\documentclass[reqno]{article}
\usepackage{amsmath}
\usepackage{amsfonts}
\usepackage{amssymb}
\usepackage{enumerate}
\usepackage[utf8]{inputenc}

\author{
Fazia {\sc Bedouhene}\footnote{Mouloud Mammeri University of Tizi-Ouzou,
Laboratoire de Math\'ematiques Pures et Appliqu\'ees, Tizi-Ouzou, Algeria
E-Mail: fbedouhene@yahoo.fr},
Nouredine {\sc Challali}\footnote{
Mouloud Mammeri University of Tizi-Ouzou,
Laboratoire de Math\'ematiques Pures et Appliqu\'ees, Tizi-Ouzou, Algeria
E-Mail: challalin@yahoo.fr},
Omar {\sc Mellah}
\footnote{Mouloud Mammeri University of Tizi-Ouzou,
Laboratoire de Math\'ematiques Pures et Appliqu\'ees, Tizi-Ouzou, Algeria
E-Mail: omellah@yahoo.fr},\\
Paul {\sc Raynaud de Fitte}
\footnote{Normandie Univ, Laboratoire Rapha\"el Salem,
UMR CNRS 6085, Rouen, France
E-Mail: prf@univ-rouen.fr},
and
Mannal {\sc Smaali}\footnote{
Mouloud Mammeri University of Tizi-Ouzou,
Laboratoire de Math\'ematiques Pures et Appliqu\'ees, Tizi-Ouzou, Algeria
E-Mail: smaali\_manel@yahoo.fr},
}

\newcommand\proof{\par\noindent{\bf Proof}\ \ }
\newcommand\proofof[1]{\par\noindent{\bf Proof of #1}\ \ }
\newcommand\finproof{{\ }\hfill\rule{2mm}{2mm}\allowbreak\medskip\par}

\newtheorem{theorem}{Theorem}[section]
\newtheorem{proposition}[theorem]{Proposition}

\newtheorem{lemma}[theorem]{Lemma}
\newtheorem{Remark}[theorem]{Remark}
\newenvironment{remark}{\begin{Remark}\em}{\end{Remark}}
\newtheorem{Example}[theorem]{Example}
\newenvironment{example}{\begin{Example}\em}{\end{Example}}
\newtheorem{Question}{Question}

\newcommand\N{\mathbb{N}}
\newcommand\R{\mathbb{R}}
\newcommand\Q{\mathbb{Q}}

\renewcommand\epsilon{\varepsilon}

\newcommand\tq{;\,} %% "such that"
\newcommand\un[1]{\,\rlap{{1}}\kern.22em \mbox{l}_{#1}} %% fonction
\newcommand\CCO[1]{\left( #1 \right)}
\newcommand\norm[1]{\left\Vert #1 \right\Vert}
\newcommand\abs[1]{\left\vert #1 \right\vert}
\newcommand\accol[1]{\left\{#1\right\}}

\newcommand\CUB{\mbox{\rm CUB}}

\newcommand\esprob{\Omega}
\newcommand\tribu{\mathcal{F}}
\newcommand\prob{\mathop{\text{\rm P}}\nolimits}

\newcommand\expect{\mathop{\text{\rm E}}\nolimits}

\newcommand\dist{\mathfrak{d}} %%{\mbox{\rm dist}} % distance generique
\newcommand\distt{\underline{\mathfrak{d}}} % une distance sur \Cont_k
\newcommand\diam{\mathop{\mbox{\textrm{Diam}}}}

\newcommand\BL{\mbox{\rm BL}}
\newcommand\bl{\mbox{\tiny\rm BL}}
\newcommand\lr{\mbox{\tiny\rm L}}

\newcommand\transl[1]{\Tilde{#1}}

\newcommand\image[1]{\mathop{{#1}_{\sharp}}}

\newcommand\Cov{\mathop{\text{\rm Cov}}\nolimits} % covariance
\newcommand\Var{\mathop{\text{\rm Var}}\nolimits} % variance

\newcommand\loi{\text{\rm law}}%{{\mathcal{L}}}
\newcommand\law[1]{ \loi\CCO{#1} } % loi de #
\newcommand\laws[1]{{\mathcal M}^{1,+}\CCO{{#1}}}

\newcommand\ellp[1]{\mathop{\text{\rm L}}\nolimits^{#1}}

\newcommand\espX{{\mathbb X}}
\newcommand\espY{{\mathbb Y}}

\newcommand\h{\mathbb{H}}

\newcommand\Cont{\text{\rm C}} % espace des fonctions continues
\newcommand\BCont{\text{\rm BC}} % espace des fonctions bornées continues
\newcommand\Dom{\mathop{\mbox{\rm Dom}}} %domaine

\newcommand\limaa[1]{\widehat{#1}} % fonction limite dans la def de aa

\newcommand\AlA{\mbox{\textrm{AA}}} % fonctions presque automorphes
\newcommand\AlAu{\mbox{\textrm{AA}}_c} % fonctions p. automorphes veech-fink
\newcommand\StAlA[2][]{\mbox{\rm{AA}}_{\St^{#2}_{#1}}} %fonctions presque
\newcommand\ER{\mathcal{E}}     % fonctions ergodiques
\newcommand\StER[2][]{\mathcal{E}_{\St^{#2}_{#1}}}
\newcommand\PAA{\mbox{\textrm{PAA}}}% fonctions pseudo presque automorphes
\newcommand\StPAA[2][]{\mbox{\textrm{PAA}}_{\St^{#2}_{#1}}}

\newcommand\AP{\mbox{\textrm{AP}}} % fonctions presque periodiques

\newcommand\AlAUb{\mbox{\textrm{AAU}}_b} % fonctions presque automorphes
\newcommand\AlAUc{\mbox{\textrm{AAU}}_c} % fonctions presque automorphes
\newcommand\ERUb{\mathcal{E}\mbox{\textrm{U}}_b}  % fonctions ergodiques U
\newcommand\ERUc{\mathcal{E}\mbox{\textrm{U}}_c}  % fonctions ergodiques U
\newcommand\PAAUb{\mbox{\textrm{PAAU}}_b}% fonctions pseudo presque
\newcommand\PAAUc{\mbox{\textrm{PAAU}}_c}
\newcommand\AAD{\mbox{\textrm{AAD}}} % fonctions pr. automorphes en distr.
\newcommand\AADm{\mbox{\textrm{AAD}}_1}
\newcommand\AADf{\mbox{\textrm{AAD}}_f}

\newcommand\PAAD{\mbox{\textrm{PAAD}}}% fonctions pseudo pr. autom. en distr.
\newcommand\PAADm{\mbox{\textrm{PAAD}}_1}
\newcommand\PAADf{\mbox{\textrm{PAAD}}_f}

\newcommand\AAcUc{\mbox{\textrm{AA}}_c\mbox{\textrm{U}}_c} %fonctions
\newcommand\AAcUb{\mbox{\textrm{AA}}_c\mbox{\textrm{U}}_b} %fonctions
\newcommand\St{\mathbb{S}} % Stepanov
\newcommand\We{\mathbb{W}} % Weyl
\newcommand\Bes{\mathbb{B}} % Besicovitch

\newcommand\UC{\textbf{C}}

\newcommand\TTm{\mbox{\rm TT}$_1$} % pte de Tudor & Tudor 1-dimensionnelle

\newcommand\trace{\mathop{\mbox{\rm tr}}}

\newcommand\Mp{\mathfrak{M}_p} % une borne dans L^p
\newcommand\fg{\mathfrak{E}} % un majorant..
\newcommand\fm{\mathfrak{f}}
\newcommand\gm{\mathfrak{g}}
\newcommand\hm{\mathfrak{h}}

\newcommand\Ymes{\mathfrak{m}} % une mesure limite dans la def de Tudor
\allowdisplaybreaks % permet de changer de page à l'intérieur des
\begin{document}
\numberwithin{equation}{section}

% \author{Fazia BEDOUHENE}
% \address{
% Department of Mathematics
% Faculty of Sciences
% University of Tizi-Ouzou, Algeria}
% \email{fbedouhene@yahoo.fr}

% \author{Nordine CHALLALI}
% \address{
% Department of Mathematics
% Faculty of Sciences
% University of Tizi-Ouzou, Algeria}
% \email{challalin@yahoo.fr}

% \author{Omar MELLAH}

% \address{Department of Mathematics
% Faculty of Sciences
% University of Tizi-Ouzou, Algeria
% }
% \email{omellah@yahoo.fr}

% \author{Paul RAYNAUD de FITTE}
% \address{Normandie Univ, Laboratoire Rapha\"el Salem,
% UMR CNRS 6085, Rouen, France}
% \email{prf@univ-rouen.fr}

% \author{Mannal SMAALI}
% \address{
% Department of Mathematics
% Faculty of Sciences
% University of Tizi-Ouzou, Algeria}
% \email{smaali\_manel@yahoo.fr}

\title{Almost automorphy and
various extensions
for stochastic processes}

% \subjclass[2010]{60G05 (60B05 60B11 60H10 34C27)}
% \keywords{Square-mean almost automorphic;
% pseudo almost automorphic in quadratic mean;
% pseudo almost automorphic in distribution;
% Stepanov;
% Ornstein-Uhlenbeck;
% semilinear stochastic differential equation;
% stochastic evolution equation
% }

\maketitle
%%%%%%%%%%%%%%%%%%%%%%%%%%%%%%%%%%%%%%%%%%%%%%%%%%%%%%%%%%%%%%%%%%%

\begin{abstract}
We compare different modes of pseudo almost automorphy and variants for stochastic processes: in probability, in
quadratic mean, or in distribution in various senses. We show by a counterexample that square-mean (pseudo)
almost automorphy is a property which is too strong for stochastic differential equations (SDEs). Finally, we
consider two semilinear SDEs, one with almost automorphic coefficients and the second one with pseudo almost
automorphic coefficients, and we prove the existence and uniqueness of a mild solution which is almost
automorphic in distribution in the first case, and pseudo almost automorphic in distribution in the second case.
\end{abstract}



%  Primary Classification
% 60 Probability theory and stochastic processes [For additional applications, see 11Kxx, 62-XX, 90-XX, 91-XX, 92-XX, 93-XX, 94-XX]
% 60G Stochastic processes
% 60G05 Foundations of stochastic processes

% Secondary Classification
% 60 Probability theory and stochastic processes [For additional applications, see 11Kxx, 62-XX, 90-XX, 91-XX, 92-XX, 93-XX, 94-XX]
% 60B Probability theory on algebraic and topological structures
% 60B05 Probability measures on topological spaces

% 60 Probability theory and stochastic processes [For additional applications, see 11Kxx, 62-XX, 90-XX, 91-XX, 92-XX, 93-XX, 94-XX]
% 60B Probability theory on algebraic and topological structures
% 60B11 Probability theory on linear topological spaces [See also 28C20]

% 60 Probability theory and stochastic processes [For additional applications, see 11Kxx, 62-XX, 90-XX, 91-XX, 92-XX, 93-XX, 94-XX]
% 60H Stochastic analysis [See also 58J65]
% 60H10 Stochastic ordinary differential equations [See also 34F05]

% 34 Ordinary differential equations
% 34C Qualitative theory [See also 37-XX]
% 34C27 Almost and pseudo-almost periodic solutions

%\tableofcontents

\section{Introduction}
Almost automorphic functions, introduced by Bochner,
are an important
generalization of almost periodic functions.
Almost automorphy is a property of regularity and recurrence of functions,
which has been studied in
the context of differential equations and dynamical systems,
and in other contexts.
The question of studying
the concept of almost automorphic
stochastic processes arises naturally in connection with stochastic
differential equations.

As can be seen in Tudor's survey \cite{Tudor95ap_processes},
almost periodicity forks into many different notions when applied
to stochastic processes: almost periodicity in probability, in
$p$-mean, in one-dimensional distributions, in finite dimensional
distributions, in distribution, almost periodicity of moments, etc.
These notions are not all comparable: for example, almost periodicity
in distribution does not imply almost periodicity in probability, and
the converse implication is false too \cite{bedouhene-mellah-prf2012}.
The situation for almost automorphy is no different.

Furthermore, things become even more complicated if one takes into account
different generalizations of almost automorphy (which have
their analogue in almost periodicity):
on the one hand, changing the mode of convergence in the definition
leads to the notions of
Stepanov-like, Weyl-like and Besicovitch-like almost automorphy.
On the other hand, for the study of asymptotic properties of
functions, it is natural to consider functions which are the sum of an
almost automorphic function and of a function vanishing in some sense
at infinity. In this way, one gets the notions of
asymptotically almost automorphic functions
and of
pseudo almost automorphic functions and their weighted variants.

Each of these notions can be interpreted in different manners for
stochastic processes, exactly in the same way as for ``plain'' almost
periodicity and almost automorphy.
It is an
objective of this paper to clarify the hierarchy of those various
concepts.
We did not try to list all possible variants, one can imagine many
other extensions and combinations, as the reader will probably do.
This is rather a preliminary groundwork, in which we investigate some
notions we think particularly useful.
For example, when we describe the different modes of
pseudo almost automorphy in distribution,
we concentrate on a stronger notion, 
which is not purely ``distributional''
but seems to be the relevant one for stochastic differential
equations.

As we have in view applications in spaces of probability
measures, which do not have a vector space structure and whose
topology can be described by different non-uniformly equivalent
metrics,
we are especially interested in the properties of almost automorphy
and pseudo almost automorphy which are purely topological, i.e.~which
do not depend on a vector structure or on a particular metric, but
only on the topology of the underlying space.

A natural application of these concepts is the study of stochastic
differential equations with almost automorphic or more general
coefficients. We provide two examples of  stochastic semilinear
evolution equations, with almost automorphic coefficients for the
first one, and
with pseudo almost automorphic coefficients for the second one,
whose unique bounded solution is
almost automorphic (respectively pseudo almost automorphic) in
distribution.
It is another objective of this paper to point out a common error in
many papers which claim the existence of nontrivial solutions which are
almost automorphic in quadratic mean. We show by a counterexample
borrowed from \cite{MRF13} that this claim is false, even for several
extensions of almost automorphy.

\paragraph{Historical comments}
%%%%%%%%
These comments are not intended to provide
a full historical account, only to
highlight some steps in the history of almost automorphic stochastic
processes and their generalizations.

The study of almost periodic random functions has a relatively long
story, starting from Slutsky \cite{Slutsky38} in 1938,
who focused on conditions for weakly stationary random processes
to have almost periodic trajectories in Besicovitch's sense.

 %\cite{bochner56stationarity}
 The investigation of almost periodicity in probability
was initiated later by the Romanian
 school \cite{Onicescu-Istratescu75,cenusa-sacuiu80,Precupanu82}.

It is only in the late eighties that almost periodicity in
distribution was considered, again by the Romanian school
(mainly by Constantin Tudor), in
connection with the study of
stochastic diferential equations with almost periodic coefficients
\cite{halanay87,halanay-tudor-morozan87,%
Morozan-Tudor89,Tudor92affine,Arnold-Tudor98,DaPrato-Tudor95}.

Starting from 2007, many papers appeared, claiming the existence of
square-mean almost periodic solutions to almost periodic semilinear stochastic
evolution equations, using a fixed point method. Despite the
counterexamples given in \cite{omarPHD,MRF13}, new papers in this vein
continue to be published.

The story of almost automorphy and its generalizations is much shorter.
Almost automorphic functions were invented by Bochner since 1955
\cite{bochner55betti,bochner61uniform,bochner62new_approach}
(the terminology stems from the fact that they were first encountered
in \cite{bochner55betti}
in the context of differential geometry on real or complex manifolds).
Almost automorphic stochastic processes and their generalizations
seem to have been investigated only
since 2010, starting with \cite{fu-liu2010square-mean-aa},
which was followed by many other papers.
Most of these papers claim almost automorphy
(or one of its generalizations)
in square mean for solutions to stochastic equations.
There are only few papers we are aware of \cite{fu2012,fu-chen13,liu-sun2014}
which investigate almost automorphy in a distributional sense.

Recently, the notion of almost
automorphic random functions in probability has been introduced by Ding,
Deng and N'Gu\'er\'ekata
\cite{ding-deng-nguerekata14}.

\paragraph{Organization of the paper}
In Section \ref{sec:generaldef}, we
present the concept of almost automorphy and some of
its generalizations:
$\mu$-pseudo almost automorphy,
Stepanov-like, Weyl-like and Besicovitch-like $\mu$-pseudo almost automorphy.
Our setting is that
of functions of a real variable with values in a %%complete
metrizable space. Metrizability seems a sufficiently general frame to
investigate almost automorphy in many useful spaces of probability
measures, while avoiding complications.
An extension to uniformizable spaces would be useful for applications
in locally convex vector spaces, this could be done using projective
limits of metrizable spaces as in \cite{bedouhene-mellah-prf2012}.
We show that almost automorphy
and a slightly generalized notion of pseudo almost automorphy
can be defined in a topological way,
without any reference to a metric nor to a vector structure.

In Section \ref{sec:stochproc}, we
investigate several notions of
almost automorphy and pseudo almost automorphy for stochastic
processes.
First, we investigate almost automorphy and its variants in $p$th mean:
the stochastic processes are seen as almost automorphic
(or, more generally, $\mu$-pseudo almost automorphic)
functions from $\R$ to $\ellp{p}$, $p\geq
0$ ($p=0$ corresponds to almost automorphy and its variants in
probability).
We show with the simple counterexample of Ornstein-Uhlenbeck process
that even a one-dimensional linear equation with constant
coefficients has no nontrivial solution which is almost automorphic
(in any of the variants considered) in $p$th mean.

Then we move to almost automorphy in distribution and its variants.
There are at least three kinds of almost
automorphy in distribution:
in one-dimensional distributions,
in finite dimensional distributions,
and in distribution of the whole process.
In the deterministic case, the first two notions are equivalent to
almost automorphy, the third one is equivalent to compact automorphy,
a stronger notion. 
For $p>0$, 
we introduce also 
the notion of almost automorphy in $p$-distribution,
which is obtained by 
adding to the preceding notions a condition of $p$-uniform
integrability.   
% Taking into account the $p$-integrability of the process ($p> 0$),
% these three notions give rise to three notions of 
%%
For $\mu$-pseudo almost automorphy,
the situation becomes even more complicated, because there are several ways to
take into account the ergodic part. We introduce the notion of
processes which are $\mu$-pseudo almost
automorphic in $p$-distribution, which are the sum of a process which
is almost automorphic in
$p$-distribution and a process which is $\mu$-ergodic in $p$th mean.  We
use this notion in the next section.
We do not address the notions of Stepanov-like, Weyl-like or Besicovitch-like (pseudo) almost automorphy for
stochastic processes, these notions would probably have to be linked to a particular choice of a metric on a
space of probability measures.

We study the superposition operator (also called Nemytskii operator)
between spaces of processes which are almost
periodic (compact almost automorphic, and $\mu$-pseudo compact almost
automorphic) in distribution.

We also carry out 
a comparison of the main notions of (generalized) almost automorphy for
stochastic processes: in
probability, in $p$th mean, and in $p$-distribution.

Finally, in Section \ref{sec:SDE}, we consider two semilinear
stochastic evolution equations in a Hilbert space.
The first one has almost automorphic coefficients, and the second one
has $\mu$-pseudo almost automorphic coefficients. 
We show that each equation has a unique mild solution
which is almost automorphic in $2$-distribution in the first case, and
$\mu$-pseudo almost automorphic in $2$-distribution in the second case.

\section{Weighted pseudo almost automorphy
in Banach spaces and in
  metric spaces}\label{sec:generaldef}
%%%%%%%%%%%%%%%%%%%%%
\subsection{Notations and definitions }

%%%%%%%%%%%%%%%%%%%%%%%%%%%%%%%%%%%%%%%%%%%%%%%%%%%%%%%%%%%%%%%
In the sequel, $\espX$ and $\espY$ are %complete
metrizable topological spaces.
When no confusion may arise,
we denote by $\dist$ a distance on $\espX$ (respectively on $\espY$)
 which generates the topology of $\espX$ (respectively $\espY$).
Most of our results depend only on the topology of those
 spaces, not on the choice of particular metrics.
% the distance on $\espX$ as well as that on
% $\espY$ is denoted by $\dist$.
When $\espX$ and $\espY$ are Banach spaces,
their norms are indistinctly denoted by $\norm{.}$, and
$\dist$ is assumed to result
from $\norm{.}$.

We denote by $\Cont(\espX,\espY)$ the space of continuous functions
from $\espX$ to $\espY$.
When this space is endowed with the topology of uniform convergence on
compact subsets of $\espX$, it is denoted by $\Cont_k(\espX,\espY)$.

For a continuous function $f :\,\R\rightarrow\espX$,
we define its {\em translation mapping}
$$\transl{f} :\,
\left\{\begin{array}{lcl}
\R&\rightarrow&\Cont(\R,\espX)\\
t&\rightarrow&f(t+.).
\end{array}\right.
$$

\subsection{Almost periodicity and almost automorphy}
\paragraph{Almost periodicity}
We say that a continuous function
$f :\,\R\rightarrow\espX$
is {\em almost periodic} if, for any $\epsilon>0$,
 there exists $l(\epsilon)>0$  such that any
 interval of length $l(\epsilon)$ contains at least
 an {\em$\epsilon$-almost period}, that is,
 a number $\tau$
 for which
$$
\dist \CCO{f(t+\tau),f(t)}\leq\varepsilon, \text{ for all } t\in\R.
$$
We denote by $\AP(\R,\espX)$ the space of $\espX$-valued
almost periodic functions.

By a result of Bochner \cite{bochner62new_approach},
$f :\R\rightarrow\espX$ is almost periodic
if, and only if, the set
$\accol{\transl{f}(t) , \,t\in \R} =
\accol{f(t+.),\,t\in \R}$
is totally bounded in the space $\Cont(\R,\espX)$ endowed with the
norm $\norm{.}_\infty$ of uniform convergence.

Another very useful characterization
is {\em Bochner's double sequence criterion} \cite{bochner62new_approach}:
$f$ is almost periodic
if, and only if, it is continuous and,
for every pair of sequences $(t'_n)$ and $(s'_n)$ in
  $\R$, there are subsequences $(t_n)$ of $(t'_n)$ and
  $(s_n)$ of $ (s'_n)$ respectively, with same indexes, such
  that, for every $t\in \R$, the limits
\begin{equation}\label{def:Bochner double sequence}
 \lim_{n\rightarrow\infty}\lim_{m\rightarrow\infty}f(t+t_n+s_m)
 \text{ and }\lim_{n\rightarrow\infty}f(t+t_n+s_n),
\end{equation}
exist and are equal.
This very useful criterion shows that the set $\AP(\R,\espX)$ depends
only on the topology of $\espX$, i.e.~it does not
depend on any uniform structure on $\espX$,
in particular it does not depend on the choice of any norm
(if $\espX$ is a vector space) or any distance on $\espX$.

\paragraph{Almost automorphy}
Almost automorphic functions were introduced by Bochner
\cite{bochner62new_approach} and studied in depth by Veech \cite{veech65aa},
see also the monographs
\cite{shen-yi98semiflows,nguerekata01book,nguerekata05book}
for applications to differential equations.
A continuous mapping $f :\,\R\rightarrow\espX$
is said to be {\em almost automorphic}
if, for every
sequence $(t'_n)$ in $\R$, there exists a subsequence $(t_n)$
such that, for every $t\in\R$, the limit
\begin{equation}\label{eq:aa1}
g(t)=\lim_{n\rightarrow\infty} f(t+t_n)
\end{equation}
exists and
\begin{equation}\label{eq:aa2}
\lim_{n\rightarrow\infty}g(t-t_n)=f(t).
\end{equation}
The range $R_f$ of $f$ is then relatively compact, because we can
extract from every
sequence $(f(t_n))$ in $R_f$ a convergent subsequence.

 Clearly, the space of almost automorphic $\espX$-valued functions
depends only on the
topology of $\espX$.

Almost automorphic functions generalize almost periodic functions in
the sense that $f$ is almost periodic
if, and only if, the above limits are uniform with respect to $t$.

Note that, in \eqref{eq:aa1} and \eqref{eq:aa2}, the limit function
$g$ is not necessarily continuous.
Let us consider the following property:
\begin{itemize}
\item[(\UC)] For any choice of $(t'_n)$ and $(t_n)$, the function
$g$ of \eqref{eq:aa1} and \eqref{eq:aa2} is continuous.
\end{itemize}
Functions satisfying (\UC) are called
{continuous almost automorphic} functions in \cite{veech65aa}.
They were re-introduced by Fink \cite{fink68aa-ap} under the
name of {\em compact almost automorphic} functions. This terminology is now
generally adopted, so we stick to it.

It has been shown by Veech in
\cite[Lemma 4.1.1]{veech65aa}
(see also \cite[Theorem 2.6]{nguerekata05comments}
and \cite{nguerekata-pankov08stepanov,lizama-nguerekata10})
that, if $f$ satifies (\UC), $f$ is
uniformly continuous.
The proof of Veech is given in the case when $\espX$ is the field of
complex numbers,
but it extends to any metric space.
Furthermore,
$f$ satifies (\UC) if, and only if,
the convergence in
\eqref{eq:aa1} and \eqref{eq:aa2} is uniform on the compact
intervals.
We denote by $\AlAu(\R,\espX)$ the subspace of functions satisfying
(\UC).

We have the inclusions
$$
\AP(\R,\espX)\subset\AlAu(\R,\espX)\subset\AlA(\R,\espX).
$$
All these spaces depend only on the topological structure of
$\espX$ and not on its metric.

\paragraph{Almost automorphic functions depending on a parameter}
Following \cite{liang-zhang_xiao-jun08},
we say that a function $f :\,\R\times\espY\rightarrow \espX$ is {\em almost automorphic with respect to the
first variable, uniformly with respect to the second variable in bounded subsets of $\espY$} (respectively {\em
in compact subsets of $\espY$}) if, for every sequence $(t'_n)$ in $\R$, there exists a subsequence $(t_n)$ such
that, for every $t\in\R$ and every $y\in\espY$, the limit
$$g(t,y)=\lim_{n\rightarrow\infty} f(t+t_n,y)$$
exists and,
for every  bounded (respectively compact) subset $B$ of $\espY$,
the convergence is uniform with respect to $y\in B$, and if the convergence
$$\lim_{n\rightarrow\infty}g(t-t_n,y)=f(t,y)$$
holds uniformly with respect to $y\in B$.
We denote by $\AlAUb(\R\times\espY,\espX)$ and
$\AlAUc(\R\times\espY,\espX) $ respectively the spaces of
such functions.

Similarly, one can define the spaces of functions
$f :\,\R\times\espY\rightarrow \espX$
which are {\em compact almost automorphic with respect to the first
  variable, uniformly
with respect to the second variable in bounded} (or {\em in compact})
{\em subsets of $\espY$}.
We denote these spaces by $\AAcUc(\R\times \espY,\espX)$ and
$\AAcUb(\R\times \espY,\espX)$ respectively.

These notions are different from the notion of functions almost automorphic
uniformly in $y$
defined in \cite{blot-al09superposition,blot-cieutat_ezzinbi2012}.

%%%%%%%%%%%%%%%%%%
\begin{proposition}\label{prop:AAcUc}
Let $f\in\AAcUc(\R\times \espY,\espX)$. Assume that $f$ is continuous
with respect to the second variable. Then $f$ is continous on $\R\times \espY$,
and, for every compact subset $K$ of $\espY$,
$f$ is uniformly continuous on $\R\times K$.
\end{proposition}
%%%%%%%%%%%%%%%%%
\proof
For simplicity, we use the
same notation $\dist$ for distances on $\espY$ and $\espX$ which
generate the topologies of $\espY$ and $\espX$ respectively.

{\em First step } Let us show that $f$ is jointly continuous. Let $(t,x)\in \R\times\espY$, and let $(t_n,x_n)$
be a sequence in $\R\times\espY$ which converges to $(t,x)$. Let $\epsilon>0$. The set $K=\{x_n\tq
n\in\N\}\cup\{x\}$ is compact, thus there exists $N_1\in\N$ such that, for any $y\in K$,
$$n,m\geq N_1\Rightarrow \dist(f(t_n,y),f(t_m,y))<\epsilon/3.$$
Now, there exists $N_2\in\N$ such that
$$n\geq N_2\Rightarrow
\dist(f(t_{N_1},x_n),f(t_{N_1},x))<\epsilon/3.$$
We deduce, for $n\geq (N_1\vee N_2)$,
\begin{multline*}
\dist(f(t_n,x_n),f(t,x))
\leq
\dist(f(t_n,x_n),f(t_{N_1},x_n))\\
+\dist(f(t_{N_1},x_n),f(t_{N_1},x))
+\dist(f(t_{N_1},x),f(t,x))
<\epsilon/3+\epsilon/3+\epsilon/3=\epsilon,
\end{multline*}
which proves the continuity of $f$.

\medskip
{\em Second step } Let $(t'_n)$ be a sequence in $\R$.
Let $(t_n)$ be a subsequence of $(t'_n)$ such that, for every
$y\in\espY$, and for every $t\in\R$, the limit
\begin{equation*}
g(t,y)=\lim_{n\rightarrow\infty} f(t+t_n,y)
\end{equation*}
exists, uniformly with respect to $y$ in compact subsets of $\espY$, and
\begin{equation*}%\label{eq:aac2}
\lim_{n\rightarrow\infty}g(t-t_n,y)=f(t,y).
\end{equation*}
By our hypothesis, for each $y\in\espY$, the function $g(.,y)$ is continuous. A similar reasoning to that of the
first step shows that $g$ is continuous on $\R\times\espY$. Indeed, let $(t,x)\in \R\times\espY$, and let
$(s_k,x_k)$ be a sequence in $\R\times\espY$ such that $(t+s_k,x_k)$ converges to $(t,x)$. Let $\epsilon>0$. The
set $K=\{x_n\tq n\in\N\}\cup\{x\}$ is compact, thus there exists an integer $N$ such that, for every $y\in K$,
$$n\geq N\Rightarrow \dist\Bigl(g(t,y),f(t+t_n,y)\Bigr)<\epsilon/3.$$
By continuity of $f$ at the point $(t+t_N,x)$, there exists $N'\in\N$ such that
$$k\geq N'\Rightarrow \dist\Bigl(f(t+s_k+t_N,x_k),f(t+t_N,x)\Bigr)<\epsilon/3.$$
We have thus, for $k\geq N'$,
\begin{align*}
\dist\Bigl(g(t+s_k,x_k),g(t,x)\Bigr)
\leq & \dist\Bigl(g(t+s_k,x_k),f(t+s_k+t_N,x_k)\Bigr)\\
&+\dist\Bigl(f(t+s_k+t_N,x_k),f(t+t_N,x)\Bigr)\\
&+\dist\Bigl(f(t+t_N,x),g(t,x)\Bigr)\\
<&\epsilon/3+\epsilon/3+\epsilon/3=\epsilon.
\end{align*}

\medskip
{\em Third step } Let $K$ be a compact subset of $\espY$.
Assume that $f$ is not uniformly continuous on $\R\times K$.
We can find two sequences $(s_n,x_n)$ and $(t_n,y_n)$ in $\R\times
K$ such that $(s_n-t_n)+\dist(x_n,y_n)$
converges to $0$ and $\dist\bigl(f(s_n,x_n),f(t_n,y_n)\bigr)>2\delta$ for
some $\delta>0$ and for all $n\in\N$.
By compactnes of $K$, and  extracting if necessary a subsequence,
we can assume that $(x_n)$ and $(y_n)$ converge to a common limit $x\in K$.
We have thus
\begin{multline*}
\liminf_{n\rightarrow\infty}\dist\Bigl(f(s_n,x_n),f(s_n,x)\Bigr)
+\liminf_{n\rightarrow\infty}\dist\Bigl(f(t_n,y_n),f(t_n,x)\Bigr)\\
\geq
\liminf_{n\rightarrow\infty}\dist\Bigl(f(s_n,x_n),f(t_n,y_n)\Bigr)>2\delta,
\end{multline*}
which implies that at least one term in the left hand side is greater
than $\delta$.
So, we can assume, without loss of generality, that
\begin{equation}\label{eq:ddelta}
\liminf_{n\rightarrow\infty}\dist\Bigl(f(t_n,y_n),f(t_n,x)\Bigr)>\delta.
\end{equation}
Extracting if necessary a further subsequence, we can assume also that
there exists a function $g\,:\R\times \espY\rightarrow\espX$ such
that
$$\lim_{n\rightarrow\infty} f(t_n,y)=g(0,y)$$
uniformly with respect to $y\in K$. We have proved in the second step
that $g$ is continuous.
But, then, we have
\begin{multline*}
\limsup_{n\rightarrow\infty}\dist\Bigl(f(t_n,y_n),f(t_n,x)\Bigr)\\
\begin{aligned}
\leq & \limsup_{n\rightarrow\infty}\Bigl\lgroup
       \dist\Bigl(f(t_n,y_n),g(0,y_n)\Bigr)\\
      & \phantom{\limsup}+\dist\Bigl(g(0,y_n),g(0,x)\Bigr)
       +\dist\Bigl(g(0,x),f(t_n,x)\Bigr)
                              \Bigr\rgroup
=0,
\end{aligned}
\end{multline*}
which contradicts \eqref{eq:ddelta}.
\finproof
%%%%%%%%%

\subsection{(Weighted) pseudo almost automorphy}
Pseudo almost periodic functions were invented by Zhang
\cite{zhang94integr,zhang94,zhang95,zhang01book}.
The generalization of this concept to pseudo almost automorphic
functions was investigated in
\cite{liang-zhang_xiao-jun08}.
To define pseudo almost automorphy, we need another class of functions.
Assume for the moment that $\espX$ is a Banach space. Let
\begin{equation*}
\ER(\R,\espX)=\accol{f\in\BCont(\R,\espX)\tq
\lim_{r\rightarrow\infty}\frac{1}{2r}\int_{[-r,r]}\norm{f(t)}\,dt=0},
\end{equation*}
where $\BCont(\R,\espX)$ denotes the space of bounded continuous
functions from $\R$ to $\espX$.
We say that a continuous function $f :\,\R\rightarrow\espX$ is
{\em pseudo almost automorphic} if it has the form
\begin{equation}
f=g+\Phi, \quad g\in\AlA(\R,\espX), \quad \Phi\in\ER(\R,\espX).
\label{eq:paa}
\end{equation}
The space of $\espX$-valued pseudo almost automorphic functions is
denoted by $\PAA(\R,\espX)$.

Weighted pseudo almost automorphic functions were introduced
by Blot et al.~in
\cite{blot-mophou-nguerekata-pennequin2009} and later generalized in
\cite{blot-cieutat_ezzinbi2012}.
They generalize the weighted pseudo almost periodic functions introduced
by Diagana \cite{diagana06weighted,diagana07book,diagana08weighted},
see also \cite{blot-cieutat-ezzinbi2013}.
Let $\mu$ be a Borel measure on $\R$ such that
\begin{equation}\label{eq:mu}
\mu(\R)=\infty \text{ and }
\mu(I)<\infty \text{ for every bounded interval $I$}.
\end{equation}
We define the space $\ER(\R,\espX,\mu)$
of {\em$\mu$-ergodic} $\espX$-valued functions
by
\begin{equation*}
\ER(\R,\espX,\mu)=\accol{f\in\BCont(\R,\espX)\tq
\lim_{r\rightarrow\infty}\frac{1}{\mu([-r,r])}\int_{[-r,r]}\norm{f(t)}\,d\mu(t)=0}.
\end{equation*}
The space $\PAA(\R,\espX,\mu)$ of
{\em $\mu$-pseudo almost automorphic functions}
with values in $\espX$
is the space of continuous functions $f :\,\R\rightarrow\espX$ of the
form
\begin{equation}
f=g+\Phi, \quad g\in\AlA(\R,\espX), \quad \Phi\in\ER(\R,\espX,\mu).
\label{eq:wpaa}
\end{equation}
The space $\PAA(\R,\espX,\mu)$ contains the
{\em asymptotically almost automorphic} functions, that is,
the functions of the form
\begin{equation*}
f=g+\Phi, \quad g\in\AlA(\R,\espX), \quad
\lim_{\abs{t}\rightarrow\infty}\norm{\Phi(t)}=0,
%\label{eq:apaa}
\end{equation*}
see \cite[Corollary 2.16]{blot-cieutat_ezzinbi2012}.

Note that, contrarily to \eqref{eq:paa},
the decomposition \eqref{eq:wpaa}
is not necessarily
unique \cite[Remark 4.4 and Theorem 4.7]{blot-cieutat_ezzinbi2012},
even in the case of weighted almost periodic functions
\cite{liang-xiao-zhang10decomposition,blot-cieutat-ezzinbi2013}.
A sufficient condition of uniqueness of the decomposition is
that $\ER(\R,\espX,\mu)$ be {translation invariant}.
This is the case in particular if
Condition (\textbf{H}) of \cite{blot-cieutat_ezzinbi2012} is satisfied:
\begin{itemize}
\item[(\textbf{H})] For every $\tau\in\R$, there exist $\beta>0$ and a
  bounded interval $I$ such that $\mu(A+\tau)\leq\beta \mu(A)$
  whenever $A$ is a Borel subset of $\R$ such that $A\cap I=\emptyset$.
\end{itemize}

The following elementary lemma will prove useful.
\begin{lemma}\label{lem:subsequence_ergodic}
Let $f\in\ER(\R,\espX,\mu)$, with $\mu$ satisfying \eqref{eq:mu}. There exists a sequence $(t_n)$ in $\R$ such
that $(\abs{t_n})$ converges to $\infty$ and $(f(t_n))$ converges to 0. If furthermore
\begin{equation}\label{eq:R+}
\liminf_{r\rightarrow\infty}\frac{\mu([0,r])}{\mu([-r,r])}>0,
\end{equation}
one can choose $(t_n)$ converging
to $+\infty$.
\end{lemma}
%%%%%%%%%%%
\proof
Observe first that, for every bounded interval $I$,
the function $f$ satisfies
\begin{equation*}
\lim_{r\rightarrow\infty}\frac{1}{\mu([-r,r]\setminus I)}
\int_{[-r,r]\setminus I}\norm{f(t)}\,d\mu(t)=0
\end{equation*}
(see \cite[Theorem 2.14]{blot-cieutat_ezzinbi2012} for a stronger result).
Assume that the first part of the lemma is false.
There exist $\epsilon>0$ and $R>0$
such that, for $\abs{t}\geq R$, $\norm{f(t)}\geq \epsilon$.
Then we have
\begin{multline*}
0=\lim_{r\rightarrow\infty}\frac{1}{\mu([-r,r]\setminus [-R,R])}
\int_{[-r,r]\setminus [-R,R]}\norm{f(t)}\,d\mu(t)\\
\geq  \lim_{r\rightarrow\infty}
\epsilon\frac{\mu([-r,r]\setminus [-R,R])}{\mu([-r,r]\setminus
  [-R,R])}
=\epsilon,
\end{multline*}
a contradiction.

In the case when \eqref{eq:R+} is satisfied,
we have
\begin{multline*}
\limsup_{r\rightarrow\infty}\frac{1}{\mu([0,r]\setminus [0,R])}
\int_{[0,r]\setminus [0,R]}\norm{f(t)}\,d\mu(t)\\
\begin{aligned}
\leq  &\limsup_{r\rightarrow\infty} \frac{\mu([-r,r]\setminus [-R,R])}{\mu([0,r]\setminus [0,R])} \,
 \frac{1}{\mu([-r,r]\setminus [-R,R])}
\int_{[-r,r]\setminus [-R,R]}\norm{f(t)}\,d\mu(t) \\
=&\limsup_{r\rightarrow\infty} \frac{\mu([-r,r])}{\mu([0,r])} \,
 \frac{1}{\mu([-r,r]\setminus [-R,R])}
\int_{[-r,r]\setminus [-R,R]}\norm{f(t)}\,d\mu(t) \leq 0.
\end{aligned}
\end{multline*}
Then we only need to reproduce the reasoning of the first part of the
lemma, replacing $[-r,r]$ by $[0,r]$.
\finproof
%%%%%%%%%

\paragraph{Pseudo almost automorphic functions depending on a parameter}
Let $\mu$ be a Borel measure on $\R$ satisfying \eqref{eq:mu}.
We say that a continuous function $f :\,\R\times\espY\rightarrow \espX$
is {\em $\mu$-ergodic with respect to the first variable, uniformly
with respect to the second variable in bounded subsets of $\espY$}
(respectively {\em in compact subsets of $\espY$}) if, for every $x\in\espY$,
$f(.,x)$ is $\mu$-ergodic,
and the convergence of $1/\mu{([-r,r])}\int^{r}_{-r} \norm{f(t,x)}\,d\mu(t)$ is
uniform with respect to $x$ in bounded subsets of $\espY$
(respectively  compact subsets of $\espY$).
The space
of such functions is denoted by $\ERUb(\R\times\espY,\espX,\mu)$ (respectively
$\ERUc(\R\times\espY,\espX,\mu)$).

\begin{remark}
\label{rem:ERUc}If
for each $x\in \espY$, $f(.,x) \in \ER(\R, \espX, \mu)$ and
$f(t,x)$ is continuous with respect to $x$, uniformly with respect to $t$,
then $f\in\ERUc(\R\times\espY,\espX,\mu)$.

Indeed, let $K$ be a compact subset of $\espY$, and let
$\epsilon>0$. Let $\dist$ be any distance on $\espY$ which generates
the topology of $\espY$.
There exists $\eta>0$ such that, for all $x,y\in K$ satisfying
$\dist\CCO{x,y}<\eta$,
we have $\norm{f(t,x)-f(t,y)}<\epsilon$ for every $t\in\R$.
Let $x_1,\dots,x_m$ be a finite sequence in $K$ such that
$K\subset \cup_{i=1}^{m}B(x_i,\eta)$. We have, for every $r>0$,
\begin{multline*}
\sup_{x\in K}\frac{1}{\mu{([-r,r])}}\int^{r}_{-r}
  \norm{f(t,x)}\,d\mu(t)\\
\begin{aligned}
\leq &\max_{1\leq i\leq m}\sup_{x\in B(x_i,\eta)}\biggl\lgroup
       \frac{1}{\mu{([-r,r])}}\int^{r}_{-r}
           \norm{f(t,x)-f(t,x_i)}\,d\mu(t)\\
      &\phantom{\max_{1\leq i\leq m}\sup_{x\in B(x_i,\eta)}\biggl\lgroup}
       +\frac{1}{\mu{([-r,r])}}\int^{r}_{-r}\norm{f(t,x_i)}\,d\mu(t)\biggr\rgroup\\
\leq &\epsilon + \max_{1\leq i\leq m} \frac{1}{\mu{([-r,r])}}\int^{r}_{-r}\norm{f(t,x_i)}\,d\mu(t),
\end{aligned}
\end{multline*}
which shows that $f\in\ERUc(\R\times\espY,\espX,\mu)$.
\end{remark}

We say that a continuous function $f :\,\R\times\espY\rightarrow\espX$ is
{\em $\mu$-pseudo almost automorphic with respect to the first variable, uniformly
with respect to the second variable in bounded subsets of $\espY$}
(respectively {\em in compact subsets of $\espY$})
if it has the form
\begin{align*}
&f=g+\Phi, \quad g\in\AlAUb(\R\times\espY,\espX),
\quad \Phi\in\ERUb(\R\times\espY,\espX,\mu)\\
%\label{eq:paau}
 \text{(respectively }
&f=g+\Phi, \quad g\in\AlAUc(\R\times\espY,\espX),
\quad \Phi\in\ERUc(\R\times\espY,\espX,\mu)
 \text{).}
\end{align*}
The space of such functions is
denoted by $\PAAUb(\R\times\espY,\espX,\mu)$
(respectively $\PAAUc(\R\times\espY,\espX,\mu)$).

\subsection{Stepanov, Weyl and Besicovitch-like
pseudo almost automorphy}
\paragraph{Stepanov-like pseudo almost automorphy and variants}
The notion of Stepanov-like almost automorphy was proposed by Casarino
in \cite{casarino00}. Then
Stepanov-like pseudo almost automorphy was first studied by Diagana
\cite{diagana09stepanov}.
Stepanov-like weighted pseudo almost automorphy seems to have been
investigated first and simultaneously in
\cite{xia-fan12} and
\cite{zhang-chang-nguerekata12}.

Let $p>0$.
We say that a locally $p$-integrable function
$f :\,\R\rightarrow\espX$ is
{\em $\St^p$-almost automorphic}, or
{\em Stepanov-like almost automorphic}
if, for every
sequence $(t_n)$ in $\R$, there exists a subsequence $(t'_n)$
and a locally $p$-integrable function $g :\,\R\rightarrow\espX$
such that, for all $t\in\R$,
$$\lim_{n\rightarrow\infty} \abs{f(.+t'_n)-g}_{\St^p}(t)=0$$
and
$$\lim_{n\rightarrow\infty}\abs{g(.-t'_n)-f}_{\St^p}(t)=0,$$
where, for any locally $p$-integrable function $h :\,\R\rightarrow\espX$,
\begin{equation*}
\abs{h}_{\St^p}(t)=
                 \CCO{\int_{t}^{t+1}\norm{h(s)}^p\,ds}^{1/p},
\end{equation*}
and assuming also that $\norm{f}_{\St^p}:=
\sup\{\abs{f}_{\St^p}(t)\tq t\in\R\}<\infty$. 
The space of $\St^p$-almost automorphic $\espX$-valued functions
is denoted by $\StAlA{p}(\R,\espX)$.

The {\em Bochner transform}\footnote{The terminology is due to the fact
  that Bochner was the first to use this transform,
  for Stepanov almost periodicity, in
\cite{bochner33abstrakte}.}
of a function $f :\,\R\rightarrow\espX$ is the
function
$$f^b :\,
\left\{ \begin{array}{lcl}
\R&\rightarrow&\espX^{[0,1]}\\
t&\mapsto & f(t+.).
\end{array}\right.
$$
We have
$$\StAlA{p}(\R,\espX)
=\accol{f \tq f^b\in\AlA(\R,\ellp{p}([0,1],dt,\espX))}.$$
We define the Stepanov-like $\mu$-ergodic functions
in a similar way:
$$\StER{p}(\R,\espX,\mu)
=\accol{f \tq f^b\in\ER(\R,\ellp{p}([0,1],dt,\espX),\mu)}.$$

Let $\mu$ be a Borel measure on $\R$ satisfying \eqref{eq:mu}.
We say that $f :\,\R\rightarrow\espX$ is
{\em $\St^p$-pseudo almost automorphic}, or
{\em Stepanov-like weighted pseudo almost automorphic}
if $f$ has the form
\begin{equation*}
f=g+\Phi, \quad g\in\StAlA{p}(\R,\espX), \quad \Phi\in\StER{p}(\R,\espX,\mu).
\end{equation*}

A further extension has been imagined by Diagana
\cite{diagana12generalized_stepanov}:
it consists in adding a weight in the Stepanov seminorms
$\abs{.}_{\St^p}$ and $\norm{.}_{\St^p}$.
Let $p$ and $\mu$ as before,
and let $\nu$ be a Borel measure on the interval $[0,1]$ such that
\begin{equation}\label{eq:nu}
0<\nu([0,1])<+\infty.
\end{equation}
Set
\begin{equation*}
\abs{h}_{\St^p_\nu}(t)=
                 \CCO{\int_{0}^{1}\norm{h(t+s)}^p\,d\nu(s)}^{1/p},
\end{equation*}
and $\norm{h}_{\St^p_\nu}=\sup\{\abs{h}_{\St^p_\nu}(t)\tq t\in\R\}$.
Then, by replacing 
$\abs{.}_{\mathbb{S}^p}$ by $\abs{.}_{\mathbb
{S}^p_\nu}$ and 
$\norm{.}_{\St^p}$ by $\norm{.}_{\St^p_\nu}$,
 one defines in the obvious
way the space
$\StAlA[\nu]{p}(\R,\espX)$
of
{\em $\St^p_\nu$-almost automorphic} $\espX$-valued functions,
the space
$\StER[\nu]{p}(\R,\ellp{p}(\esprob,\prob,\R),\mu)$
of {\em $\St^p_\nu$-$\mu$-ergodic functions},
and the space $\StPAA[\nu]{p}(\R,\espX)$
of
{\em $\St^p_\nu$-pseudo almost automorphic} %$\espX$-valued
functions.

\paragraph{Weyl-like and Besicovitch-like pseudo almost automorphy}
The concept of Weyl-like pseudo almost automorphy
has been recently explored by Abbas
\cite{abbas14weyl}. The definition is similar to that of Stepanov-like
pseudo almost automorphy, replacing $\abs{.}_{\St^p}$ and
$\norm{.}_{\St^p}$ 
by the
weaker seminorms
$\abs{.}_{\We^p}$ and 
$\norm{.}_{\We^p}$, defined by
\begin{align*}
&\abs{h}_{\We^p}(x)=\lim_{r\rightarrow+\infty}
                 \CCO{\frac{1}{2r}\int_{x-r}^{x+r}\norm{h(t)}^p\,dt}^{1/p}\\
\intertext{and}
&\norm{h}_{\We^p}=\lim_{r\rightarrow+\infty}\sup_{x\in\R}
                 \CCO{\frac{1}{2r}\int_{x-r}^{x+r}\norm{h(t)}^p\,dt}^{1/p}.
\end{align*}
A further weakening leads to the Besicovitch seminorm, which does not seem to
have been investigated in the context of almost automorphy:
\begin{equation*}
\norm{h}_{\Bes^p}=\limsup_{ r \rightarrow+\infty}
                 \CCO{\frac{1}{2r}\int_{-r}^{r}\norm{h(t)}^p\,dt}^{1/p}
\end{equation*}
(here, $\abs{h}_{\Bes^p}(x)=\norm{h}_{\Bes^p}$ 
for all $x\in\R$).
We shall briefly consider this seminorm in Example \ref{exple:OU}.

%%%%%%%%%%%%%%%%%%%%%%%%%%%%%%%%%%%%%%%%%%%%%%%%%%%%%%%%%%%%
\subsection{Weighted pseudo almost automorphy
in topological spaces}
We have seen that, to define the space of almost automorphic
functions,
as well as that of almost periodic functions, on a space $\espX$,
there is no need to assume that $\espX$ is a vector space, nor a
metric space,
these spaces depend only on the topological structure of $\espX$.
We prove in this section that
the definition of $\PAA(\R,\espX,\mu)$ is metric (independent of the
vector structure of $\espX$ but dependent on the metric),
and that it can be made purely topological if one allows for a
slight change of the definition.
This leads to two topological concepts of $\mu$-pseudo
  almost automorphy: in C.~and M.~Tudor's sense, and in the wide sense.

In this section, unless otherwise stated, $\espX$ is only assumed to be a
  topological space, not necessarily metrizable.

%%%%%%%%%%%%%%%%%%%%%%%%%%%%%%%%%%%%%%%%
\begin{remark}\label{rem:ER_topological}
Assume that $\espX$ is a Banach space.
By definition, the space $\ER(\R,\espX,\mu)$ consists of
continuous functions $f :\,\R\rightarrow\espX$ satisfying
\begin{enumerate}[(i)]
\item\label{cond:fbounded} $f$ is bounded,
\item\label{cond:ergodic}
$\displaystyle \lim_{r\rightarrow\infty}\frac{1}{\mu([-r,r])}\int_{[-r,r]}\norm{f(t)}\,d\mu(t)=0.$
\end{enumerate}
Condition \eqref{cond:fbounded} is of metric nature.
By \cite[Theorem 2.14]{blot-cieutat_ezzinbi2012},
Condition \eqref{cond:ergodic} implies Condition
\eqref{eq:ER_topological} below, and the converse implication is true
if we assume \eqref{cond:fbounded} :
\begin{equation}\label{eq:ER_topological}
\text{For any $\epsilon>0$, }
\lim_{r\rightarrow\infty}
  \frac{\mu{\{t\in[-r,r]\tq \norm{f(t)}>\epsilon\}}}{\mu([-r,r])}=0.
\end{equation}
Thus $f$ satisfies \eqref{eq:ER_topological} if, and only if,
\begin{equation*}%\label{eq:ergobound}
\lim_{r\rightarrow\infty}\frac{1}{\mu([-r,r])}
           \int_{[-r,r]}
 \bigl\lgroup\norm{f(t)}\wedge 1\bigr\rgroup\,d\mu(t)=0.
\end{equation*}
Condition \eqref{eq:ER_topological} can be reformulated as
\begin{equation}\label{eq:ER_topological-neighbourhood}
\text{For any neighbourhood $U$ of $0$, }
\lim_{r\rightarrow\infty}
  \frac{\mu{\{t\in[-r,r]\tq f(t)\not\in U\}}}{\mu([-r,r])}=0.
\end{equation}
The vector structure of $\espX$ is still involved in
\eqref{eq:ER_topological-neighbourhood} through the vector
$0$. But we use $\ER(\R,\espX,\mu)$ only in order to ensure that
a function is close in a certain sense to  $\AlA(\R,\espX)$.
To that end,
as $\AlA(\R,\espX)$ contains the constant functions,
we can allow a generalization of $\ER(\R,\espX,\mu)$ by
replacing $0$ in
\eqref{eq:ER_topological-neighbourhood} by any other fixed point $x_0$
of $\espX$.\footnote{
Actually, the terminology ``ergodic'' is misleading but it has the advantage of shortness. It would be more
appropriate to follow Zhang's terminology (e.g.~\cite{zhang94,zhang01book}), and call {\em $\mu$-ergodic
perturbations} the elements of $\ER(\R,\espX,\mu)$, and call {\em ergodic} the functions $f$ such that
$1/{\mu([-r,r])}\int^r_{-r} f(t)\,d\mu(t)$ converges to some limit, not necessarily 0. }
\end{remark}
%%%%%%%%%%%%

An elegant metric definition of
$\mu$-pseudo almost periodic functions
has been proposed by Constantin and Maria Tudor in \cite{Tudor-Tudor99}.
If we adapt their definition,
a continuous function $f :\,\R\rightarrow\espX$ is
{\em $\mu$-pseudo almost automorphic in C.~and M.~Tudor's sense}
if $f$ has relatively compact
range and there exists $g\in\AlA(\R,\espX)$ such that
the function $t\mapsto \dist(f(t),g(t))$ is in $\ER(\R,\R,\mu)$.
This definition is more restrictive than the standard one
because
the metric condition \eqref{cond:fbounded} is replaced
by the stronger topological condition that
$f$ have relatively compact range
(note that a subset $A$ of $\espX$ is relatively compact if, and only
if, it is bounded for {\em every} metric which generates the topology
of $\espX$,
see \cite[Problem 4.3.E.(c)]{Engelking89}).

Let us say that a continuous function $f :\,\R\rightarrow\espX$ is
{\em $\mu$-pseudo almost automorphic in the wide sense}
if there exists $g\in\AlA(\R,\espX)$ such that
the function $t\mapsto \dist(f(t),g(t))\wedge 1$ is in $\ER(\R,\R,\mu)$.
We can get rid of the distance $\dist$ in this definition
by using the fact that the closure of the range
of any function in $\AlA(\R,\espX)$ is compact.
On a compact space $K$, there is one and only one uniform structure,
a basis of entourages of which
consists of all sets of the form
$V=\cup_{i=1}^m \CCO{U_i\times U_i}$, where $U_1,\dots,U_m$ is a
finite open cover of $K$,
see e.g.~\cite{Bourbaki71-I-IV} or \cite{Engelking89}.
In this way, we obtain \eqref{cond:paa-topol} below,
which shows that the space of functions which are
$\mu$-pseudo almost automorphic in the wide sense
depends only on the topology of $\espX$.

We have thus two possible topological definitions of $\mu$-pseudo
almost automorphy: in C.~and M.~Tudor's sense, and in the wide
sense. The former is stronger than \eqref{eq:wpaa},
while the latter is weaker.

%%%%%%%%%%%%%%%%%%%%%%%%%%%%%%%%%%%%%%%%%
\begin{proposition}\label{prop:equivalence}
(Topological characterization of $\mu$-pseudo almost automorphy in the
wide sense)
Let $f : \R\rightarrow\espX$
be continuous. Let $\mu$ be a Borel measure on $\R$ satisfying \eqref{eq:mu}.
If $\espX$ is a metric space, the following propositions are equivalent.
\begin{enumerate}[(i)]

\item\label{cond:paa-metric}
$f$ is $\mu$-pseudo almost automorphic in the wide
  sense,
i.e.~there exists a function $g\in\AlA(\R,\espX)$ such that
$$\lim_{r\rightarrow\infty}\frac{1}{\mu([-r,r])}
           \int_{[-r,r]}
 \bigl\lgroup\dist\CCO{f(t),g(t)}\wedge 1\bigr\rgroup\,d\mu(t)=0.$$
\item\label{cond:paa-topol}
There exists a function $g\in\AlA(\R,\espX)$ such that,
for any finite open cover $U_1,\dots,U_m$
of the closure $K$ of ${\{ g(t) \tq t\in\R \}}$,
$$\lim_{r\rightarrow\infty}
  \frac{\mu{\{t\in[-r,r]\tq \CCO{f(t),g(t)}\not\in
      V\}}}{\mu([-r,r])}=0,$$
where $V=\cup_{i=1}^m \CCO{U_i\times U_i}$.
\end{enumerate}

\end{proposition}
%%%%%%
\proof

$\eqref{cond:paa-metric} \Rightarrow \eqref{cond:paa-topol}$. Recall that $K$ is compact because
$g\in\AlA(\R,\espX)$. Let $V=\cup_{i=1}^m \CCO{U_i\times U_i}$ be as in \eqref{cond:paa-topol}.
Then $V$ is an open neighborhood in $\espX\times\espX$ of the diagonal $\Delta=\{(x,x)\tq x\in K\}$. Define the
distance $\dist_2$ on $\espX\times\espX$ by
$$\dist_2\Bigl((x,y),(x',y')\Bigr)=\dist(x,x')+\dist(y,y').$$
As $\Delta$ is compact, the distance $\epsilon=\dist_2(\Delta,\espX \times \espX\setminus V)$ is positive. Let
$$B_\epsilon:=\{(y,z)\in\espX\times\espX
         \tq \dist_2\Bigl((y,z),\Delta\Bigr) <\epsilon\}\subset V.$$
For each $(y,z)\in B_\epsilon$, there exists $x\in K$ such that
$\dist_2\Bigl((y,z),(x,x)\Bigr)<\epsilon$, thus
$$\dist(y,z)\leq \dist(y,x)+\dist(x,z)<\epsilon.$$
On the other hand, if $z\in K$ and $\dist(y,z)<\epsilon$, we have
$$\dist_2\Bigl((y,z),(z,z)\Bigr)=\dist(y,z)<\epsilon,$$
thus $(y,z)\in B_\epsilon$.
Applying this to $(f(t),g(t))$, we get,
for every $r>0$,
\begin{equation*}
  \frac{\mu{\{t\in[-r,r]\tq \CCO{f(t),g(t)}\not\in
      V\}}}{\mu([-r,r])}
\leq
  \frac{\mu{\{t\in[-r,r]\tq
           \dist\CCO{f(t),g(t)}\geq\epsilon\}}}{\mu([-r,r])}.
\end{equation*}
But \eqref{cond:paa-metric} means that the function
$t\mapsto\dist(f(t),g(t))\wedge 1$ is in $\ER(\R,\R,\mu)$, thus,
by \cite[Theorem 2.14]{blot-cieutat_ezzinbi2012}, the latter term goes
to 0 when $r$ goes to $\infty$.

\medskip
$\eqref{cond:paa-topol} \Rightarrow \eqref{cond:paa-metric}$.
Let $\epsilon>0$.
Let $U_1,\dots,U_m$ be a finite open cover of $K$ such that
$\diam(U_i)>\epsilon$, $i=1,\dots,m$,
and let $V=\cup_{i=1}^m\CCO{U_i\times U_i}$.
We have
\begin{equation*}
  \frac{\mu{\{t\in[-r,r]\tq
      \dist\CCO{f(t),g(t)}>\epsilon\}}}{\mu([-r,r])}
\leq
  \frac{\mu{\{t\in[-r,r]\tq \CCO{f(t),g(t)}\not\in
      V\}}}{\mu([-r,r])},
\end{equation*}
where the latter term goes to $0$ when $r$ goes to $\infty$.
The conclusion follows from \cite[Theorem 2.14]{blot-cieutat_ezzinbi2012}.
\finproof
%%%%%%%%%

%%%%%%%%%%%%%%%%%%%%%%%%%%%%%%%%%%%%%
\begin{remark}\label{rem:class-topol}
The reasoning of Proposition \ref{prop:equivalence} can be applied without
change to give a topological characterization of
$\mu$-pseudo almost periodic functions in the wide sense.
More generally,
$\AlA(\R,\espX)$ can be replaced in this reasoning
by any class of functions which have relatively compact range.
\end{remark}

%%%%%%%%%%%%%%% THEO UNIQUENESS %%%%%%%%%%%%%
\begin{theorem}\label{theo:uniqueness}
(Uniqueness of the decompostion of pseudo almost automorphic functions
in the wide sense)
Let  $\mu$ be a Borel measure  on $\R$
satisfying \eqref{eq:mu} and Condition (\textbf{H}).
Let $f :\,\R\rightarrow\espX$ satisfying
Condition \eqref{cond:paa-topol} of Proposition \ref{prop:equivalence}.
Then the function $g\in\AlA(\R,\espX)$
given by Condition \eqref{cond:paa-topol} is unique and satisfies
$$\{g(t)\tq t\in\R\}\subset \overline{\{f(t)\tq t\in\R\}}.$$
\end{theorem}
%%%%%%%%%%%%%
\proof
Let $g$ and $K$ as in Condition \eqref{cond:paa-topol}.
Let $\varphi :\,\espX\rightarrow\R$ be a continuous function.
Then $\varphi\circ g\in \AlA(\R,\R)$.
Let $\Hat{U}_1,\dots,\Hat{U}_m$ be a finite open cover
of the closure $\varphi(K)$ of ${\{ \varphi\circ g(t) \tq t\in\R
  \}}$, and let $\Hat{V}=\cup_{i=1}^m \CCO{\Hat{U}_i\times \Hat{U}_i}$.
Then $U_1=\varphi^{-1}(\Hat{U}_1),\dots,U_m=\varphi^{-1}(\Hat{U}_m)$
form a finite open cover of $K$.
Let $V=\cup_{i=1}^m \CCO{U_i\times U_i}$. We have
\begin{align*}
  \frac{\mu{\{t\in[-r,r]\tq \CCO{\varphi\circ f(t),\varphi\circ g(t)}\not\in
      \Hat{V}\}}}{\mu([-r,r])}
\leq &
  \frac{\mu{\{t\in[-r,r]\tq \CCO{f(t),g(t)}\not\in
      V\}}}{\mu([-r,r])}.
\end{align*}
Thus, by Proposition \ref{prop:equivalence},
$\varphi\circ f$ is in $\PAA(\R,\R,\mu)$ and the function
$\varphi\circ f-\varphi\circ g$ is in $\ER(\R,\R,\mu)$.
By \cite[Theorem 4.1]{blot-cieutat_ezzinbi2012},
we have
$$\{\varphi\circ g(t)\tq t\in\R\}
\subset \overline{\{\varphi\circ f(t)\tq t\in\R\}},$$
and,
by the uniqueness of the decomposition of vector-valued $\mu$-pseudo
almost automorphic functions
\cite[Theorem 4.7]{blot-cieutat_ezzinbi2012} % uniquenessse
if $g'\in\AlA(\R,\espX)$ satisfies the same condition as $g$,
we have $\varphi\circ g'=\varphi\circ g$.
As $\varphi$ is arbitrary, we deduce $g'=g$.

\finproof
%%%%%%%%%

%%%%%%%%%%%%%%%%%%%%%%%%%%%%%%%%%%%%%%%%%%%%%%%%%%%%%%%%%%%%%
\section{Weighted pseudo almost automorphy
for stochastic processes}
\label{sec:stochproc}

From now on, $\espX$ and $\espY$ are assumed to Polish spaces,
i.e.~separable metrizable topological spaces whose topology is
generated by a complete metric.

%%%%%%%%%%%%%%%%%%%%%%%%%%%%%%%%%%%%%%%%%%%%%%%%%%%%%%%%%%%%%%%%%%%
\subsection{Weighted pseudo almost automorphy in $p$th mean}
We assume here that $\espX$ is a Banach space.
Let $X=(X_t)_{t\in\R}$ be a continuous stochastic process with
values in $\espX$, defined on a probability space $(\esprob,\tribu,\prob)$.
Let $\mu$ be a Borel measure on $\R$ satisfying \eqref{eq:mu}.

Let $p>0$. We say that $X$ is {\em almost automorphic in $p$th mean}
(respectively {\em
$\mu$-pseudo almost automorphic in $p$th mean}) if the mapping
$t\mapsto X(t)$ is in $\AlA(\R,\ellp{p}(\esprob,\prob,\espX))$
(respectively in $\PAA(\R,\ellp{p}(\esprob,\prob,\espX),\mu)$,
i.e., if it has the form $X=Y+Z$,
where $Y\in \AlA(\R,\ellp{p}(\esprob,\prob,\espX))$
and $Z \in \ER(\R,\ellp{p}(\esprob,\prob,\espX),\mu)$).
When $p=2$, we say that $X$ is {\em square-mean almost
automorphic} (respectively {\em square-mean $\mu$-pseudo almost automorphic}).

The process $X$ is said to be {\em almost automorphic in probability}
if the mapping
 $X:\,t \rightarrow \ellp{0}(\esprob,\prob,\espX)$ is almost
 automorphic,
where $\ellp{0}(\esprob,\prob,\espX)$ is the space of measurable
mappings from $\esprob$ to $\espX$, endowed with the topology of
convergence in probability.
Recall that the topology of $\ellp{0}(\esprob,\prob,\espX)$
is induced by e.g.~the distance
$$\dist_{\text{\tiny Prob}}(U,V)=\expect \CCO{\norm{U-V}\wedge1},$$
which is complete.

The process $X$ is said to be
{\em $\mu$-pseudo almost automorphic in probability},
and we write $X \in
\PAA(\R,\ellp{0}(\esprob,\prob,\espX),\mu)$,
if the mapping $t\mapsto X(t)$,
               $\R\rightarrow \ellp{0}(\esprob,\prob,\espX)$
is $\mu$-pseudo almost automorphic in the wide sense (or, equivalently, if it is $\mu$-pseudo almost automorphic
when $\ellp{0}(\esprob,\prob,\espX)$ is endowed with $\dist_{\text{\tiny Prob}}$), i.e.~if it has the form
$X=Y+Z$ where $Y \in\AlA(\R,\ellp{0}(\esprob,\prob,\espX))$ and
$Z$ %$Z\in\ER$
satisfies
\begin{equation}\label{eq:ergodic0}
\lim_{r\rightarrow\infty}\frac{1}{\mu([-r,r])}
           \int_{[-r,r]} \expect
 \bigl\lgroup\norm{Z(t)}\wedge 1\bigr\rgroup\,d\mu(t)=0.
\end{equation}
We denote by
$\ER(\R,\ellp{0}(\esprob,\prob,\espX),\mu)$ the set
of stochastic processes $Z$ satisfying \eqref{eq:ergodic0}.

Note that, for $p\geq 0$,  the previous decompositions of
$X\in\PAA(\R,\ellp{p}(\esprob,\prob,\espX),\mu)$  are unique
under the condition \textbf{(H)}.

Clearly, for $0\leq p\leq q$,
we have
$$\AlA(\R,\ellp{q}(\esprob,\prob,\espX))\subset
\AlA(\R,\ellp{p}(\esprob,\prob,\espX))$$ and
$$\PAA(\R,\ellp{q}(\esprob,\prob,\espX),\mu)
\subset \PAA(\R,\ellp{p}(\esprob,\prob,\espX),\mu).$$
Conversely,
if the set
$\{\norm{X(t)}^q\tq t\in \R\}$ is uniformly integrable, we have the
implications
\begin{gather*}
  \bigl\lgroup X\in\AlA(\R,\ellp{p}(\esprob,\prob,\espX)) \bigr\rgroup
  \Rightarrow
  \bigl\lgroup X\in\AlA(\R,\ellp{q}(\esprob,\prob,\espX))\bigr\rgroup,\\
  \bigl\lgroup X\in\PAA(\R,\ellp{p}(\esprob,\prob,\espX),\mu) \bigr\rgroup
  \Rightarrow
  \bigl\lgroup X\in\PAA(\R,\ellp{q}(\esprob,\prob,\espX),\mu)\bigr\rgroup.
\end{gather*}

A process $X$ is
{\em Stepanov-like almost automorphic in $p$th mean}
if $X$ is in $\StAlA{p}(\R,\ellp{p}(\esprob,\prob,\espX))$.
We define in the same way the processes which are
{\em Stepanov-like $\mu$-pseudo almost automorphic in $p$th mean},
{\em Weyl-like ($\mu$-pseudo) almost automorphic in $p$th mean},
{\em Besicovitch-like ($\mu$-pseudo) almost automorphic in $p$th mean}.

%%%%%%%%%%%%%%%%%%%%%%%%%%%%%%%%%%%%%%%%%%%%%%%%%%%%%%%%%%%%%%%%%%%
\paragraph{Explicit counterexample
to square-mean pseudo almost automorphy}
At the time of the submission of this paper, there are at least 24 papers
listed in Mathematical Reviews %%%on MathSciNet
related to almost automorphy of solutions to
stochastic differential equations, which all have been published in
2010 or later.
To our knowledge, except for \cite{fu2012,fu-chen13,liu-sun2014}, all other papers
claim the existence of square-mean pseudo almost automorphic solutions
to stochastic differential equations with coefficients having similar
properties.

We show that a  very simple counterexample from
\cite{omarPHD,MRF13} contradicts these claims.
The other counterexamples given in \cite{omarPHD,MRF13}
also contradict these claims.

%%%%%%%%%%%%%%%%%%%%%%%%%%%%%%
\begin{example}\label{exple:OU}{\rm
{\bf (Stationary Ornstein-Uhlenbeck process)}
Let $W=(W(t))_{t\in\R}$ be a standard Brownian motion on the real
line.
Let $\alpha,\sigma>0$,
and let $X$ be the stationary Ornstein-Uhlenbeck process
(see \cite{lindgren2006stationary_processes}) defined by
\begin{equation}
  \label{eq:ornstein-uhlenbeck}
X(t)=\sqrt{2\alpha\sigma}\int_{-\infty}^{t} e^{-\alpha(t-s)}dW(s).
\end{equation}
Then $X$ is the only $\ellp{2}$-bounded solution of the following SDE,
which is a particular case of Equation (3.1) in
\cite{bezandry-diagana07existence}:
\begin{equation*}
 dX(t)=-\alpha X(t)\,dt+\sqrt{2\alpha} \sigma\,dW(t).
\end{equation*}

The process $X$ is Gaussian with mean $0$,
and we have, for all $t\in\R$ and $\tau\geq 0$,
\begin{equation*}
\Cov(X(t),X(t+\tau))=\sigma^2 e^{-\alpha \tau}.
\end{equation*}

%%%%%%%%%%%%%%%%%%%%%%%%%%%%%%%%%
Assume that $X$ is
square-mean $\mu$-pseudo almost automorphic,
for some Borel measure $\mu$ on $\R$
satisfying \eqref{eq:mu} and \eqref{eq:R+}.
Then we can decompose $X$ as
$$X=Y+Z, \quad
Y\in\AlA(\R,\ellp{2}(\esprob,\R)),
\quad
Z\in\ER(\R,\ellp{2}(\esprob,\R),\mu).
$$
By Lemma \ref{lem:subsequence_ergodic}, we can find an increasing sequence $(t_n)$ of real numbers which
converges to $+\infty$ such that $(Z(t_n))$ converges to 0 in $\ellp{2}(\esprob,\R)$. Then, we can extract a
sequence (still denoted by $(t_n)$ for simplicity) such that $(Y(t_n))$ converges in $\ellp{2}$ to a random
variable $\widehat{Y}$. Thus $(X(t_n))$ converges in $\ellp{2}$ to $\widehat{Y}$. Necessarily
$\widehat{Y}$ is Gaussian with law $\mathcal{N}(0,2\alpha\sigma^2)$, and $\widehat{Y}$ is
$\mathcal{G}$-measurable, where $\mathcal{G}=\sigma\CCO{X_{t_n}\tq
  n\geq 0}$. Moreover $(X(t_n),\widehat{Y})$ is Gaussian
for every $n$, and we have, for any integer $n$,
\begin{equation*}
\Cov(X({t_n}),\widehat{Y})=\lim_{m\rightarrow\infty}\Cov(X({t_n}),X({t_{n+m}}))=0,
\end{equation*}
because
$(X^2(t))_{t\in\R}$ is uniformly integrable.
This proves that $\widehat{Y}$ is independent of $X({t_n})$ for every $n$, thus
$\widehat{Y}$ is independent of $\mathcal{G}$. Thus $\widehat{Y}$ is constant, a
contradiction.
Thus \eqref{eq:ornstein-uhlenbeck} has no
square-mean $\mu$-pseudo almost automorphic
solution.

Let us show that $X$ is not
Weyl-like nor Besicovitch-like square-mean pseudo almost automorphic.
It is enough to disprove the Besicovitch sense.
Assume that $X$ is Besicovitch-like square-mean pseudo almost automorphic.
As before, using Lemma \ref{lem:subsequence_ergodic},
we can find a sequence $(t_n)$ converging to $+\infty$ and
a process $\widehat{Y}$ such
that
$$\lim_{n\rightarrow\infty}\bigl\lVert X(t_n)-\widehat{Y}\bigr\rVert_{\Bes^2}=0.$$
In particular, $(X(t_n))$ is Cauchy for $\norm{.}_{\Bes^2}$, thus,
for every $\epsilon>0$, there exists $N(\epsilon)\in\N$, such that,
for all $n,m\in\N$,
\begin{equation}\label{eq:Cauchy-Bes}
(n\geq N(\epsilon))\Rightarrow
\bigl\lVert X(t_n)-X(t_{n+m})\bigr\rVert_{\Bes^2}{}\leq \epsilon.
\end{equation}
But we have
\begin{multline*}
\bigl\lVert X(t_n)-X(t_{n+m})\bigr\rVert_{\Bes^2}^2\\
\begin{aligned}
=&\limsup_{r\rightarrow+\infty}
      \frac{1}{2r}
      \int_{-r}^r \expect\bigl\lvert X(t_n+s)-X(t_{n+m}+s)\bigr\rvert^2 ds\\
=&\limsup_{r\rightarrow+\infty}
      \frac{1}{2r}
      \int_{-r}^r \expect\bigl( X^2(t_n+s)+X^2(t_{n+m}+s)-2X(t_n+s)X(t_{n+m}+s)\bigr)
      ds\\
=&\limsup_{r\rightarrow+\infty}
      \frac{1}{2r}
      \int_{-r}^r \bigl( \sigma^2+\sigma^2-2\sigma^2e^{-\alpha(t_{n+m}-t_n)}\bigr)
      ds\\
=&2\sigma^2\bigl(1-e^{-\alpha(t_{n+m}-t_n)}\bigr).
\end{aligned}
\end{multline*}
For $m$ large, the last term is arbitrarily close to $2\sigma^2$, which
contradicts \eqref{eq:Cauchy-Bes} for $\epsilon<2\sigma^2$.

A similar calculation shows that,
for any Borel measure $\mu$ on $\R$
satisfying \eqref{eq:mu} and \eqref{eq:R+}, and for any Borel measure
$\nu$ on $[0,1]$ satisfying \eqref{eq:nu},
the process $X$ is not square-mean $\St^2_\nu$-$\mu$-pseudo almost automorphic.
}
\end{example}

%%%%%%%%%%%%%%%%%%%%%%%%%%%%%%%%%%%%%%%%%%%%%%
\subsection{Weighted pseudo almost automorphy
in distribution}

  We denote by $\law{X}$ the law (or distribution) of a random variable
  $X$. For any topological space $\espX$, we denote by
  $\laws{\espX}$ the set of Borel probability measures on $\espX$,
  endowed with the topology of narrow (or weak) convergence,
  i.e.~the coarsest topology such that the mappings
  $\mu\mapsto\mu(\varphi)$, $\laws{\espX}\rightarrow\R$
  are continuous for all bounded continuous
  $\varphi :\,\espX\rightarrow\R$.

If $\tau :\,\espX\rightarrow\espY$ is a Borel measurable mapping and $\mu$
is a Borel measure on $\espX$, we denote by $\image{\tau}\mu$ the
Borel measure on $\espY$ defined by
$$\image{\tau}\mu(B)=\mu(\tau^{-1}(B))$$
for every Borel set of $\espY$.

Let $\BCont(\espX,\R)$ denote
the space of bounded continuous functions from $\espX$ to
$\R$, which we endow with the norm
$$\norm{\varphi}_\infty=\sup_{x\in\espX}\abs{\varphi(x)}.$$
For a given distance $\dist$ on $\espX$,
and for $\varphi \in \BCont(\espX,\R)$ we define
\begin{gather*}
\norm{\varphi}_{\lr}
= \sup\Bigl\{\dfrac{\varphi(x)-\varphi(y)}{\dist(x,y)} \tq x \neq y\Bigl\}\\
\norm{\varphi}_{\bl} = \max\{\norm{\varphi}_\infty,\norm{\varphi}_{\lr}\}.
\end{gather*}
We denote
$$\BL(\espX,\R)
= \bigl\{\varphi\in \BCont(\espX,\R);\norm{\varphi}_{\bl}<\infty
\bigl\}.$$
The {\em bounded Lipschitz distance} $\dist_{\bl}$ associated with $\dist$
on $\laws{\espX}$ is
defined by
$$\dist_{\bl}(\mu, \nu) = \sup_{\substack{\varphi\in\bl(\espX,\R)\\
\norm{\varphi}_{\bl}\leq 1}}
\int_\espX \varphi \,d(\mu-\nu)  .$$
This metric generates the narrow
(or weak) topology on $\laws{\espX}$.

Let $p\geq 0$, and let $(X_n)$ be a sequence in
$\ellp{p}(\esprob,\prob,\espX)$.
We say that $(X_n)$
{\em converges in $p$-distribution}
(or simply {\em converges in distribution} if $p=0$)
to a random vector $X$ if
\begin{enumerate}[(i)]
\item the sequence $(\law{X_n})$ converges to $\law{X}$ for the narrow topology on
  $\laws{\espX}$,
\item if $p>0$, the sequence $(\norm{X_n}^p)$ is uniformly integrable.
\end{enumerate}

\paragraph{Almost automorphy in distribution}
If $X$ is continuous with values in $\espX$, we denote by $\transl{X}(t)$
the random variable $X(t+.)$ with values in $\Cont(\R,\espX)$.

We say that $X$ is {\em almost automorphic in one-dimensional
  distributions}
if the mapping $t\mapsto
\law{{X}(t)}$,
$\R\mapsto\laws{\espX}$ is almost automorphic.

\begin{remark}\label{rem:1dim_pasbon}
One-dimensional distributions of a process reflect poorly its
properties. For example,
let $X$ be the Ornstein-Uhlenbeck process of Example \ref{exple:OU},
and set $Y(t)=X(0)$, $t\in\R$.
The processes $X$ and $Y$ have
the same one-dimensional distributions
with completely different
trajectories and behaviors.
The trajectories of $Y$ are constant, whereas the covariance
$\Cov(X(t+\tau),X(t))$ converges to $0$ when $\tau$ goes to $\infty$.
\end{remark}

We say that $X$ is {\em almost automorphic in finite dimensional distributions}
if, for every finite sequence $t_1,\dots,t_m\in\R$,
the mapping $t\mapsto
\law{X(t+t_1),\dots,X(t+t_m)}$,
$\R\mapsto\laws{\espX^m}$ is almost automorphic.

We say that $X$ is {\em almost automorphic in distribution}
if the mapping $t\mapsto \law{\transl{X}(t)}$,
$\R\mapsto\laws{\Cont_k(\R,\espX)}$ is almost automorphic,
where $\Cont_k(\R,\espX)$ denotes the space $\Cont(\R,\espX)$
endowed with the
topology of uniform convergence on compact subsets.
The Ornstein-Uhlenbeck process $X$ of Example \ref{exple:OU}
is almost automorphic in distribution because the mapping
$t\mapsto \law{\transl{X}(t)}$ is constant.

\begin{remark}\label{rem:distrib_vs_deterministic}
If $X$ is a deterministic process $\R\rightarrow \espX$, we have the
equivalences
\begin{align*}
X\in\AlA(\R,\espX)
\Leftrightarrow&
\text{$X$ is almost automorphic in one-dimensional distributions}\\
\Leftrightarrow&
\text{$X$ is almost automorphic in finite dimensional distributions}\\
\intertext{and}
X\in\AlAu(\R,\espX)
\Leftrightarrow&
\text{$X$ is almost automorphic in distribution}.
\end{align*}
\end{remark}

Actually, as Remark \ref{rem:distrib_vs_deterministic} suggests,
the definition of almost automorphy in distribution
implies a stronger property for the mapping
$t\mapsto
\law{\transl{X}(t)}$.
We need first some notations.
For simplicity, we set $\Cont_k=\Cont_k(\R,\espX)$.
For every $t\in\R$, we define a continuous operator on $\Cont_k$:
$$
\tau_t :\,
\left\{\begin{array}{lcl}
\Cont_k & \rightarrow & \Cont_k\\
x & \mapsto & x(t+.)=\transl{x}(t).
\end{array}\right.
$$

%%%%%%%%%%%%%%%
\begin{proposition}\label{prop:AAuD}
If $X$ is almost automorphic in distribution, the mapping
$t\mapsto{\law{\transl{X}(t)}}$ is in
$\AlAu(\R,\laws{\Cont_k(\R,\espX)})$.
Furthermore, for any sequence $(t_n)$ in $\R$ such that,
for every $t\in\R$,
$(\law{\transl{X}(t+t_n)})$ converges to a limit
$g(t)\in\laws{\Cont_k}$,
the function $g$ satisfies, for every $t\in\R$,
the consistency relation
\begin{equation}\label{eq:imagetau}
{g}(t)=\image{\CCO{\tau_t}}g(0).
\end{equation}
\end{proposition}
%%%%%%%%%%%
\proof
Let us denote $f(t)=\law{\transl{X}(t)}$.
We have, for every $t\in\R$,
\begin{equation*}
f(t)=\law{\tau_t\circ{X}}
=\image{\CCO{\tau_t}}f(0).
\end{equation*}
Let $(t_n)$ be a sequence in $\R$ such that, for each $t\in\R$,
the sequence $(f(t+t_n))$ converges to some $g(t)\in\laws{\Cont_k}$.
Then, for every $t\in\R$,
\eqref{eq:imagetau} is satisfied
by continuity of the operator $\image{\CCO{\tau_t}}$.
To prove the continuity of $g$, let us endow $\Cont_k$ with the
distance
\begin{equation}\label{eq:distt}
\distt(x,y)=\sum_{k\geq 1}2^{-k}\sup_{-k\leq t\leq
  k}\Bigl(\dist\bigl(x(t),y(t)\bigr)\wedge 1\Bigr),
\end{equation}
 and let $\distt_{\bl}$ be the associated bounded Lipschitz distance on
 $\laws{\Cont_k}$.
Let us show that the convergence of $(f(.+t_n))$ is
uniform for $\distt_{\bl}$ on compact intervals.
Let $r\geq 1$ be an integer. For every $t\in[-r,r]$, and for all
$x,y\in\Cont_k$, we have
\begin{align*}
\distt\CCO{\tau_{t}(x),\tau_{t}(y)}
=&\sum_{k\geq 1}2^{-k}\sup_{-k\leq s\leq k}
     \Bigl(\dist\bigl(x(s+t),y(s+t)\bigr)\wedge 1\Bigr)\\
\leq&\sum_{k\geq 1}2^{-k}\sup_{-k-r\leq s\leq k+r}
     \Bigl(\dist\bigl(x(s),y(s)\bigr)\wedge 1\Bigr)\\
\leq&2^{r}\distt(x,y).
\end{align*}
Thus, for any 1-Lipschitz mapping
$\varphi :\,\Cont_k\rightarrow\R$, the mapping $\varphi\circ\tau_t$ is
$2^{r}$-Lipschitz. We deduce that, if $\norm{\varphi}_{\bl}\leq 1$,
we have $\norm{\varphi\circ\tau_t}_{\bl}\leq 1+2^{r}$.
We have thus
\begin{align*}
\distt_{\bl}\bigl(f(t+t_n),f(t+t_{n+m})\bigr)
=&\sup_{\norm{\varphi}_{\bl}\leq 1}{\expect\Bigl(
   \varphi\circ\tau_t\circ\transl{X}(t_n)-\varphi\circ\tau_t\circ\transl{X}(t_{n+m})
                                    \Bigr)}\\
\leq&(1+2^{r})\sup_{\norm{\varphi}_{\bl}\leq 1}{\expect\Bigl(
   \varphi\circ\transl{X}(t_n)-\varphi\circ\transl{X}(t_{n+m})
                                    \Bigr)}\\
=&(1+2^{r})\distt_{\bl}\bigl(f(t_n),f(t_{n+m})\bigr),
\end{align*}
which shows that $(f(.+t_n))$ is uniformly Cauchy on $[-r,r]$.
Thus $g$ is continuous,
and $f\in\AlAu(\R,\laws{\Cont_k(\R,\espX)})$.
\finproof
%%%%%%%%%

We denote by
\begin{itemize}
\item $\AADm(\R,\espX)$ the set of $\espX$-valued processes
which are almost automorphic in
one-dimensional distributions,
\item $\AADf(\R,\espX)$ the set of $\espX$-valued processes
which are almost automorphic in
finite dimensional distributions,
\item $\AAD(\R,\espX)$ the set of $\espX$-valued processes
which are almost automorphic in distribution.
\end{itemize}
(In these notations, we omit the probability
space $(\esprob,\tribu,\prob)$, as there is no ambiguity here.)
We have the inclusions
$$
\AAD(\R,\espX)\subset \AADf(\R,\espX)\subset\AADm(\R,\espX).
$$
The following result is in the line of
\cite[Theorem 2.3]{bedouhene-mellah-prf2012}.

%%%% theo AADf => AAD
\begin{theorem}\label{theo:AAFD}
Let $X$ be an $\espX$-valued stochastic process, and let $\dist$ be a
distance on $\espX$ which generates the topology of $\espX$.
Asume that $X$ satisfies the tightness condition
\begin{multline}\label{eq:Condition(2)-BMRF}
{\forall [a,b] \subset\R,}\ {\forall \varepsilon>0},\ {\forall \eta>0},\
{\exists \delta>0,}\ {\forall r\in\R,}\ \\
  \prob\bigl\{\,
   \sup_{\begin{subarray}{l}
            |t-s| <\delta\\
             t,s \in [a,b]
   \end{subarray}
} \dist\CCO{X(r+t), X(r+s)}>\eta \bigr\} <\varepsilon.
\end{multline}
Then the following
 properties are equivalent:
\begin{enumerate}[{\rm (a)}]
 \item \label{enum:AAFD} $X\in \AADf(\R,\espX)$.
  \item\label{enum:AAD}$X\in\AAD(\R,\espX)$.
\end{enumerate}
\end{theorem}
%%%%%%%%%%%%%
\proof
Clearly (\ref{enum:AAD})$\Rightarrow$(\ref{enum:AAFD}).
Assume that $X\in \AADf(\R,\espX)$.
Let  $(\gamma'_n)$ be a sequence in $\R$,
and, for $t_1, t_2, \dots, t_k, t \in \R$
define (using notations of  \cite{Tudor95ap_processes})
$$\mu_{t}^{t_1, \dots, t_k} := \law{X(t_1+t),\dots,X(t_k+t)}.$$
By a diagonal procedure we can find a subsequence
$(\gamma_n)$ of $(\gamma'_n)$  such that,
for every $k\geq1$, for all $q_1, q_2, \dots, q_k \in
\Q\bigcap\R$ (where $\Q$ is the set of rational numbers),
and for every $t\in \R$,
\begin{equation*}
\lim_n \lim_m \mu_{t+\gamma_n-\gamma_m}^{q_1, \dots, q_k}
=\mu_{t}^{q_1, \dots, q_k}.
\end{equation*}
Let $\dist_k$ be the distance on $\espX^k$ defined by
$$\dist_k\bigl((x_1,\dots,x_k),(y_1,\dots,y_k)\bigr)
=\max_{1\leq i\leq k}\dist(x_i,y_i),$$
and let $\dist_{\bl}$
the associated bounded Lipschitz distance on $\laws{\espX^k}$.
We have, for all $t_1, t_2, \dots, t_k, t \in \R$,
for all $q_1, q_2, \dots,q_k \in \Q\bigcap\R$,
and for all $n,m\in\N$,
\begin{multline*}
\dist_{\bl}\Bigl(\mu_{t+\gamma_n-\gamma_m}^{q_1, \dots, q_k},
\mu_{t+\gamma_n-\gamma_m}^{t_1, \dots, t_k}\Bigr)
= \sup_{\parallel
  f \parallel_{\bl}\leq 1}
   {\int_{\espX^k}f
d\bigl(\mu_{t+\gamma_n-\gamma_m}^{q_1, \dots, q_k}-
       \mu_{t+\gamma_n-\gamma_m}^{t_1, \dots, t_k}\bigr)}\\
\leq \max_{1\leq i\leq k} \int_{\esprob}\dist\Bigl(X(q_i+t +\gamma_n-\gamma_m), X(t_i+t
+\gamma_n-\gamma_m)\Bigr) d\prob,
\end{multline*}
so that, by \eqref{eq:Condition(2)-BMRF},
if $(q_1, \dots, q_k) \rightarrow (t_1, \dots, t_k)$,
then
$$\dist_{\bl}\Bigl(\mu_{t+\gamma_n-\gamma_m}^{q_1, \dots, q_k},
\mu_{t+\gamma_n-\gamma_m}^{t_1, \dots, t_k}\Bigr) \rightarrow0$$
uniformly with respect to $t\in\R$ and $n,m\in\N$.
By a classical result on inversion of limits, we deduce that,
for all $k\geq1$ and $t_1, \dots, t_k, t \in \R$,
\begin{equation*}
\lim_n \lim_m \mu_{t+\gamma_n-\gamma_m}^{t_1, \dots, t_k}
=\mu_{t}^{t_1, \dots, t_k}.
\end{equation*}
Therefore, to show that
\begin{equation*}
\lim_n \lim_m \law{\transl{X}(t +\gamma_n - \gamma_m)}
=\law{\transl{X}(t)},
\end{equation*}
it is enough to prove that $(\transl{X}(t))_{t\in \R}$ is
tight in $\Cont_k(\R,\espX)$.
Since $X\in \AADf(\R,\espX)$,
the family $(X(t))_{t\in\R}=(\transl{X}(t)(0))_{t\in\R}$ is tight,
by Prokhorov's theorem for relatively compact sets of probability measures on
Polish spaces.
By \eqref{eq:Condition(2)-BMRF} and the Arzel\`a-Ascoli-type
characterization of tight subsets of $\laws{\espX}$
(see e.g.~the proof of \cite[Theorem 7.3]{Billingsley99} or \cite[Theorem 4]{Whitt70}),
we conclude
that $(\transl{X}(t))_{t\in \R}$ is tight in $\Cont_k(\R,\espX)$,
which proves our claim.
\finproof

%%%%%%%% PAS UN ESPACE VECTORIEL %%%%%%
\begin{remark}\label{rem:pasvectoriel}
Assume that $\espX$ is a vector space.
The spaces $\AAD(\R,\espX)$, $\AADf(\R,\espX)$, and $\AADm(\R,\espX)$
are not vector spaces.
Indeed, let $X$ be
the Ornstein-Uhlenbeck process of Example \ref{exple:OU}.
% The process $X$ is stationary in the strong sense,
% thus it is in $\AAD(\R,\espX)$. %almost automorphic in distribution.
For each $t\in\R$, let $Y(t)=X(0)$.
The processes $X$ and $Y$ are stationary in the strong sense,
thus they are in $\AAD(\R,\espX)$.
For each $t\in\R$, the variable $Z(t)=X(t)+Y(t)$ is Gaussian centered
with variance
\begin{align*}
\Var{Z(t)}&=\expect\CCO{X^2(t)}+\expect\CCO{Y^2(t)}
                 +2\Cov\bigl({X(t),Y(t)}\bigr)\\
&=2\sigma^2+2\sigma^2\exp\CCO{-\alpha \abs{t}}\\
&\rightarrow 2\sigma^2\text{ when }\abs{t}\rightarrow\infty.
\end{align*}
Thus $\law{Z(t)}$ is the Gaussian distribution  $\mathcal{N}(0,2\sigma^2(1+\exp\CCO{-\alpha \abs{t}}))$, which
converges when $\abs{t}\rightarrow\infty$ to $\mathcal{N}(0,2\sigma^2)$. Set
$\Ymes(t)=\mathcal{N}(0,2\sigma^2)$, $t\in\R$. For each $t\in\R$, we have
\begin{gather*}
  \lim_{n\rightarrow\infty}\law{Z(t+n)}=\Ymes(t)\\
  \lim_{n\rightarrow\infty}\Ymes(t-n)=\Ymes(t)\not=\law{Z(t)}.
\end{gather*}
Thus $Z\not\in\AADm(\R,\espX)$.

This contradicts \cite[Lemma 2.3]{fu2012}.
\end{remark}

\paragraph{Almost automorphy in $p$-distribution}
A useful variant of almost automorphy in distribution takes into
account integrability of order $p$. Let $p\geq 0$. We say that a
continuous  $\espX$-valued stochastic process is
{\em almost automorph in $p$-distribution} if
\begin{enumerate}[(i)]
\item $X\in\AAD(\R,\espX)$,
\item if $p>0$, the family $(\norm{X(t)}^p)_{t\in\R}$ is uniformly integrable.
\end{enumerate}
These conditions imply that
the mapping $t\mapsto X(t)$, $\R\rightarrow\ellp{p}(\esprob,\prob,\espX)$,
is continuous.

We denote by $\AAD^p(\R,\espX)$ the set of $\espX$-valued processes
which are almost automorphic in $p$-distribution, in particular we
have $\AAD^0(\R,\espX)=\AAD(\R,\espX)$.
Similarly, for $p\geq 0$, 
one defines the sets $\AADf^p(\R,\espX)$ and $\AADm^p(\R,\espX)$ of
processes which are respectively
{\em almost automorphic in one-dimensional $p$-distribu\-tions} and
{\em almost automorphic in finite dimensional $p$-distributions}.

\paragraph{Weighted pseudo almost automorphy in distribution and variants}
As usual, we assume that $\mu$ is a Borel measure on $\R$
satisfying \eqref{eq:mu}.

C.~and M.~Tudor proposed in \cite{Tudor-Tudor99} a very natural and
elegant notion of
pseudo almost periodicity in
  (one-dimensional) distribution
that can easily be extended to weighted
$\mu$-pseudo almost automorphy:
$X$ is {\em $\mu$-pseudo almost periodic in one-dimensional distributions
in C.~and M.~Tudor's sense} if
it satisfies
\begin{itemize}
\item[(\TTm)] The mapping $t\mapsto\law{X(t)}$ is continuous with
  relatively compact range in $\laws{\espX}$,
  and
there exists an almost automorphic function
$\Ymes:\,\R\rightarrow\laws{\espX}$ such that
\begin{equation}\label{eq:TT1}
  \lim_{r\rightarrow\infty}\frac{1}{\mu([-r,r])}
            \int_{-r}^r \dist_{\bl}(\law{X(t)},\Ymes(t))d\mu(t)=0.
\end{equation}
\end{itemize}
Similar definitions are easy to write for
$\mu$-pseudo almost automorphy in finite distributions or
in distribution.

Recall that $\dist_{\bl}$ is bounded.
If we remove the condition of relatively compact range
(as in Proposition \ref{prop:equivalence}), we get three
distributional notions of {\em $\mu$-pseudo almost automorphy in the wide sense:
in one-dimensional distributions, in finite dimensional distributions,
and in distribution.}

We propose three stronger notions of $\mu$-pseudo almost automorphy in
a distributional sense, that seem to be particularly useful for
stochastic equations.

Assume that $\espX$ is a vector space.
Let $p\geq 0$.
We say that $X$ is
{\em $\mu$-pseudo almost automorphic in  $p$-distribution}
if $X$ can be written
$$X=Y+Z, \text{ where }Y\in \AAD^p(\R,\espX) \text{ and }
Z\in\ER(\R,\ellp{p}(\esprob,\prob,\espX),\mu).$$
The set of $\espX$-valued processes which are $\mu$-pseudo almost automorphic in $p$-distribution is denoted by
$\PAAD^p(\R,\espX)$. Similar definitions hold for the spaces $\PAADm^p(\R,\espX,\mu)$ and
$\PAADf^p(\R,\espX,\mu)$ of processes which are {\em $\mu$-pseudo almost automorphic in one-dimensional
$p$-distributions} and {\em in finite dimensional  $p$-distributions} respectively.

%%%%%%%%%%%%%%%%%%%%%% TT voir Remark \ref{rem:ER_topological}
\begin{remark}
The definitions we propose for $\mu$-pseudo almost automorphy in distribution, or in finite dimensional
distributions, or in one-dimensional distributions, are in a way stronger and less natural than those in the
wide distributional sense, because they involve (at least, apparently) not only the distribution of the process,
but the probability space $(\esprob,\tribu,\prob)$. For example, a random process  $X$ is in
$\PAADm^0(\R,\espX,\mu)$ if, and only if, there exists a process $Y$ defined on the same probability space such
that $\Ymes(.):=\law{Y(.)}$ is in $\AlA(\R,\laws{\espX})$ and satisfies \eqref{eq:TT1}. Note also
that, in our definition, the ergodic part is ergodic in $p$th mean. Our definitions are thus intermediate
between $\mu$-pseudo almost automorphy in a purely distributional sense and $\mu$-pseudo almost automorphy in
$p$th mean.

However, our definitions seem to be convenient for calculations, and
Theorem \ref{theo:main} shows that, for some stochastic differential
equations with $\mu$-pseudo almost automorphic coefficients, the
process $Y$ appears naturally: it is the solution of the corresponding
SDE where the coefficients are the almost automorphic parts of the
coefficients of the original SDE.
\end{remark}

Almost automorphy in distribution and some of its variants enjoy some
stability properties, as shows the following superposition lemma.
Similar
results could be proved by the same method for one-dimensional or for
finite dimensional distributions.
%%%%%%%%%%%%%%%%%%%%%%%%%%%%%%%%%
\begin{theorem}\label{theo:superpos}
{(Superposition lemma)}
Let $X$ be a continuous
$\espX$-valued stochastic process,
and let
$f :\,\R\times\espX\rightarrow\espY$ be a continuous mapping.
\begin{enumerate}
\item If $X$ is almost periodic in distribution and
$f$ is almost periodic with respect to the first variable, uniformly
with respect to the second variable in compact subsets of $\espX$,
then $f(.,{X}(.))$ is almost periodic in distribution.

\item \label{cond:superpos-aa}
If $X$ is almost automorphic in distribution and
$f$ is compact almost automorphic
with respect to the first variable, uniformly
with respect to the second variable in compact subsets of $\espX$,
then $f(.,{X}(.))$ is almost automorphic in distribution.

\item  \label{cond:supperpo-paa}
Let $\mu$ be a Borel measure on $\R$ satisfying
  \eqref{eq:mu}, and let $p\geq 0$.
Assume that $f$ is $\mu$-pseudo compact almost automorphic with
respect to the first variable,
uniformly with
respect to the second variable in compact subsets of $\espX$, i.e.
$$
f=g+h,\text{ with }g\in\AAcUc(\R\times \espX,\espY)\text{ and }
h\in \ERUc(\R\times \espX,\espY,\mu).
$$
Assume that $g$ is continuous with respect to the second
variable.
Assume furthermore that $f$ is
uniformly continuous in the second variable in compact subsets of
$\espX$, uniformly with respect to the first
variable, that is,
\begin{multline}\label{eq:ucontK}
\begin{aligned}
&\text{for every $\epsilon>0$, and for every compact subset $K$ of
  $\espX$,} \\
&\text{there exists $\eta>0$ such that, for all $x,y\in K$,}
\end{aligned}\\
\norm{x-y}\leq \eta \Rightarrow
\sup_{t\in\R}\norm{f(t,x)-f(t,y)}\leq \epsilon.
\end{multline}
If $p>0$, 
assume also 
that $f$ and $g$ satisfy the growth condition
\begin{equation}\label{eq:growth}
  \norm{f(t,x)}+\norm{g(t,x)}\leq C(1+\norm{x})
\end{equation}
for all $(t,x)\in\R\times\espX$ and for some constant $C$, and that 
$f$ is Lipschitz with respect to the second variable, 
uniformly with respect to the first one.
If $X$ is $\mu$-pseudo almost automorphic in
$p$-distribution, then
$f(.,{X}(.))$ is
$\mu$-pseudo almost automorphic in
$p$-distribution.
\end{enumerate}
\end{theorem}
%%%%%%%%%%%%%
\proof
We only prove the second and third items,
the first one can be proved in the same way as
\ref{cond:superpos-aa}, using, for example,
Bochner's double sequence criterion.

\medskip
\ref{cond:superpos-aa}.
For each $t\in\R$, for each $x\in\Cont_k(\R,\espX)$, and for each $s\in\R$,
let us denote
$$\transl{f}(t,x)(s)=f(s+t,x(s)).$$
By continuity of $f$ on $\R\times\espX$,
$\transl{f}(t,.)$ maps $\Cont_k(\R,\espX)$ to $\Cont_k(\R,\espY)$.

Let $\mathcal{K}$ be a compact subset of $\Cont_k(\R,\espX)$.
By the Arzel\`a-Ascoli Theorem
(see e.g.~\cite[Theorems 8.2.10 and 8.2.11]{Engelking89}),
this means that
$\mathcal{K}$ is closed in $\Cont_k(\R,\espX)$ and equicontinuous,
and that, for every
compact interval $I$ of $\R$,
the set $\{x(t)\tq x\in\mathcal{K},\ t\in I\}$ has compact closure in $\espX$.
Let $t\in\R$, and let $(t_n)$ be a sequence in $\R$ converging to $t$.
Let $I$ be a compact interval of $\R$, and
let $K$ be the closure of $\{x(s)\tq x\in\mathcal K,\ s\in I\}$.
We have, for any $y\in K$, and for any $s\in\R$,
$$\lim_{n}f({t_n}+s,y) =f(t+s,y)$$
where the convergence is uniform with respect to $y\in K$ and $s\in I$,
because $f$ is compact almost
automorphic uniformly with respect to the second variable in compact
subsets of $\espX$.
In particular we have,
uniformly with respect to $x\in\mathcal{K}$ and $s\in I$,
\begin{align*}
\lim_{n}\transl{f}(t_n,x)(s)
&=\lim_{n}{f}({t_n}+s,x(s)) =f(t+s,x(s)) =\transl{f}(t,x)(s),
\end{align*}
which proves that the mapping
$\transl{f}(.,x) :\,\R\rightarrow\Cont_k(\R,\espX)$, is continuous,
uniformly with respect to $x$ in compact subsets of
$\Cont_k(\R,\espX)$.

Let us check that
$\transl{f} :\, \R\times\Cont_k(\R,\espX)\rightarrow\Cont_k(\R,\espX)$
is compact almost automorphic with respect to the first variable,
uniformly with respect to the second variable in compact subsets of
$\Cont_k(\R,\espX)$.
Let $(t'_n)$ be a sequence in $\R$.
There exists a subsequence $(t_n)$
such that, for every $t\in\R$ and every $y\in\espX$,
$$\lim_{n}\lim_{m}{f(t+t_n-t_m,y)}\\
=f(t,y),$$
where the convergence is uniform with respect to $y$ in compact
subsets of $\espX$ and $t$ in compact intervals of $\R$.
For each $t\in\R$, and for each $s\in\R$, we have
\begin{align*}
\lim_{n}\lim_{m}{\transl{f}(t+t_n-t_m,x)}(s)
=&\lim_{n}\lim_{m}{f(s+t+t_n-t_m,x(s))}\\
=&f(s+t,x(s))=\transl{f}(t,x)(s),
\end{align*}
and these convergences are uniform with respect to $t$ and $s$ in
compact intervals, and with respect to $x\in\mathcal{K}$,
which proves our claim.

%%%%%%%%%%%%%%%%%%%%%%%%%%%%%%%%%%%%%%%%%%%%%%%%%%%%%%%%%%%%%%%%%%%%%
Let $\varphi :\Cont_k(\R,\espY)\rightarrow\R$ be bounded Lipschitz,
with $\norm{\varphi}_{\bl}\leq 1$, where
  $\norm{.}_{\bl}$ is taken relatively to a distance $\distt$ which generates
 the topology of $\Cont_k(\R,\espY)$,
 for example the distance defined by \eqref{eq:distt}.
Let $(t'_n)$ be a sequence in $\R$.
Let $(t_n)$ be a subsequence such that, for every $t\in\R$ and every $y\in\espX$,
\begin{align*}
\lim_{n}\lim_{m}{f(t+t_n-t_m,y)}
&=f(t,y),\\%\text{ and }
\lim_{n}\lim_{m}\law{\transl{X}(t+t_n-t_m)}
&=\law{\transl{X}(t)}.
\end{align*}
Let $\epsilon>0$.
We can find a compact subset $\mathcal{K}_\epsilon$ of
$\Cont_k(\R,\espX)$ such that, for every $t\in\R$,
\begin{equation*}
  \prob\accol{\transl{X}(t)\in\mathcal{K}_\epsilon}\geq 1-\epsilon.
\end{equation*}
Let $t\in\R$ be fixed, and, for all $n,m$,
let $\Omega_{\epsilon,n,m}$ be the measurable subset of $\Omega$ on
which $\transl{X}(t+t_n-t_m)\in\mathcal{K}_\epsilon$.
We have
 \begin{multline*}  \Bigl\vert \expect\Bigl\lgroup
      \varphi\circ\transl{f}(t+t_n-t_m,\transl{X}(t+t_n-t_m))
     -\varphi\circ\transl{f}(t, \transl{X}(t))
       \Bigr\rgroup\Bigr\vert\\
\begin{aligned}
\leq & \Bigl\vert \expect\Bigl\lgroup
      \varphi\circ\transl{f}(t+t_n-t_m,\transl{X}(t+t_n-t_m))
     -\varphi\circ\transl{f}(t, \transl{X}(t+t_n-t_m))
       \Bigr\rgroup\Bigr\vert \\
     &+ \Bigl\vert \expect\Bigl\lgroup
      \varphi\circ\transl{f}(t,\transl{X}(t+t_n-t_m))
     -\varphi\circ\transl{f}(t, \transl{X}(t))
       \Bigr\rgroup\Bigr\vert\\
\leq &  \expect\Bigl\lgroup \un{\Omega_{\epsilon,n,m}} \distt \Big(
      \transl{f}(t+t_n-t_m,\transl{X}(t+t_n-t_m))
     ,\transl{f}(t, \transl{X}(t+t_n-t_m))\Big)
       \Bigr\rgroup \\
     &+  \expect\Bigl\lgroup  \un{\Omega_{\epsilon,n,m}^c}\Bigl\vert
      \varphi\circ\transl{f}(t+t_n-t_m,\transl{X}(t+t_n-t_m))
     -\varphi\circ\transl{f}(t,  \transl{X}(t+t_n-t_m)))\Bigr\vert
       \Bigr\rgroup\\
     &+  \Bigl\vert \expect\Bigl\lgroup
      \varphi\circ\transl{f}(t,\transl{X}(t+t_n-t_m))
     -\varphi\circ\transl{f}(t, \transl{X}(t))
       \Bigr\rgroup\Bigr\vert\\
 =    &A_{n,m}+B_{n,m}+C_{n,m}.%+D_{n,m}.
\end{aligned}
\end{multline*}
We have $$A_{n,m} \leq \expect\Bigl\lgroup \un{\Omega_{\epsilon,n,m}} \sup_{x\in
\mathcal{K}_\epsilon}\distt \Big(
      \transl{f}(t+t_n-t_m,x)
     ,\transl{f}(t, x)\Big)
       \Bigr\rgroup
\leq \sup_{x\in
\mathcal{K}_\epsilon}\distt \Big(
      \transl{f}(t+t_n-t_m,x)
     ,\transl{f}(t, x)\Big)$$
thus by the almost automorphy property of $\transl{f}$,
we have $\lim_n\lim_m A_{n,m} = 0$.
Furthermore, $B_{n,m}\leq2\prob(\Omega_{\epsilon,n,m}^c)\leq 2\epsilon$
because $\norm{\varphi}_{\bl}\leq 1$.
Finally, $\lim_n\lim_m C_{n,m} = 0$ by boundedness and continuity of
$\transl{f}(t,.) :\, \Cont_k(\R,\espX)\rightarrow \Cont_k(\R,\espX)$
and the convergence in
distribution of $\transl{X}(t+t_n-t_m)$ to $\transl{X}(t)$ for each $t\in \R$.
As $\epsilon$ and $\varphi$ are arbitrary, we have proved that
$$\lim_{n}\lim_{m}
\law{\transl{f}(t+t_n-t_m,\transl{X}(t+t_n-t_m))}
=\law{\transl{f}(t,\transl{X}(t))},$$
thus the mapping
$t\mapsto\law{\transl{f}(t,\transl{X}(t))}$
is almost automorphic in distribution.

\medskip
%%%%%%%%%%% 3.
\ref{cond:supperpo-paa}. We use ideas of the proof of \cite[Theorem 5.7]{blot-cieutat_ezzinbi2012}. 
Let $(Y,Z)$
be a decomposition of $X$, namely,
\begin{equation*}
X=Y+Z,\, Y\in\AAD^p(\R,\espX),\, Z\in\ER(\R,\ellp{p}(\esprob,\prob,\espX),\mu).
\end{equation*}
The function $f(.,X(.))$
can be decomposed as
$$f(t,X(t))=g(t,Y(t))+f(t,X(t))-f(t,Y(t))+h(t,Y(t)).$$
Let $G(t)=g(t,Y(t))$ and $H(t)=f(t,X(t))-f(t,Y(t))+h(t,Y(t))$. 
By using \ref{cond:superpos-aa}.~and
\eqref{eq:growth}, we see that $t\mapsto G(t)$ is in
$\AAD^p(\R,\espY)$. 
Furthermore, by
\eqref{eq:growth}, the continuity of $f$ and the $\mu$-pseudo almost automorphy in $p$-distribution of $X$, we
have, using Vitali's theorem, that $f(.,X(.))$ is a continuous $\ellp{p}(\esprob,\prob,\espY)$-valued function.
Indeed, if $t_n\rightarrow t$, then the sequence $(\norm{X(t_n)}^p)$ is uniformly integrable, by continuity of
the mapping $t\mapsto X(t)$, $\R\rightarrow \ellp{p}(\esprob,\prob,\espX)$, and this entails that
$(f(t_n,X(t_n)))$ is uniformly integrable. To show that $f(.,X(.))$ is in $\PAAD^p(\R,\espY)$, it
 is enough to prove that
$H \in \ER(\R,\ellp{p}(\esprob,\prob,\espY),\mu)$.

Clearly $H$ is in $\Cont(\R,\ellp{p}(\esprob,\prob,\espY))$,
and bounded in $\ellp{p}(\esprob,\prob,\espX)$ (when $p>0$), by
\eqref{eq:growth}.
As $Y$ is in $\AAD(\R,\espX)$,
the family
$(\transl{Y}(t))_{t\in\R}=({Y}(t+.))_{t\in\R}$ is uniformly tight in
$\Cont_k(\R,\espX)$. For each $\epsilon>0$, there exists a compact subset
$\mathcal{K}_\epsilon$ of $\Cont_k(\R,\espX)$ such that, for every $t\in\R$,
$$\prob\accol{\transl{Y}(t)\in\mathcal{K}_\epsilon}\geq 1-\epsilon.$$
By the Arzel\`a-Ascoli Theorem
(see e.g.~\cite[Theorems 8.2.10 and 8.2.11]{Engelking89}),
this implies that, for every $\epsilon>0$, and for every compact
interval $I$ of $\R$, there exists a compact subset $K_{\epsilon,I}$
such that, for every $t\in\R$;
$$\prob\accol{(\forall s\in I)\ Y(t+s)\in K_{\epsilon,I}}\geq 1-\epsilon.$$
In particular, the family $(Y(t))_{t\in\R}$ is tight, i.e., 
denoting $K_{\epsilon}=K_{\epsilon,\{0\}}$, we have, for every
$t\in\R$, 
\begin{equation*}%\label{eq:Ytendu-thm-super-paa}
\prob\accol{Y(t)\in K_{\epsilon}}\geq 1-\epsilon.
\end{equation*}
Let $\Omega_{\epsilon,t}$ be the measurable subset of $\Omega$ on which
$Y(t)\in K_{\epsilon}$.
The function $g$ is uniformly continuous on $\R\times K_{\epsilon}$ by
Proposition \ref{prop:AAcUc}.
We deduce by \eqref{eq:ucontK} that
there exists $\eta(\epsilon)>0$ such that,
for all
$y_1,y_2\in K_{\epsilon}$,
$$
\norm{y_1-y_2}\leq \eta(\epsilon)
\Rightarrow \sup_{t\in\R}\bigl(\norm{h(t,y_1)-h(t,y_2)}\bigr) \leq \epsilon.$$
We can find a finite sequence
$(y_i)_{1\leq i\leq m}$ in $K_{\epsilon}$
such that
$$K_{\epsilon} \subset
\bigcup_{i=1}^{m}B(y_i,\eta(\epsilon)).$$
We have, for every $t\in\R$, 
\begin{multline*}
\expect\bigl(\norm{ h(t,Y(t))}\wedge 1\bigr)\\
\begin{aligned}
\leq&\expect\biggl(\min_{1\leq i\leq m}\bigl(
      \un{\Omega_{\epsilon,t}}  \norm{ h(t,Y(t))-h(t,y_i) }\wedge 1\bigr)\biggr)
     +\max_{1\leq i\leq m}\norm{h(t,y_i)}\\
 &+\expect\bigl(\un{\Omega_{\epsilon,t}^c}\norm{h(t,Y(t))}\wedge 1\bigr)\\
\leq&\epsilon+\max_{1\leq i\leq m}\norm{h(t,y_i)}\wedge 1
    +\prob\CCO{\Omega_{\epsilon,t}^c}\\
\leq&\max_{1\leq i\leq m}\norm{h(t,y_i)} +2\epsilon.
\end{aligned}
\end{multline*}
Since for all $i\in \{1,\dots, m\}$,
the function $t\rightarrow h(t,y_i)$  satisfies
$$\lim_{r\rightarrow\infty}\frac{1}{\mu([-r,r])}
       \int_{[-r,r]}\norm{h(t,y_i)}d\mu(t)=0,$$
we deduce that, for every $\epsilon>0$,
$$
\limsup_{r\rightarrow\infty}\frac{1}{\mu([-r,r])}\int_{[-r,r]}
      \expect\bigl(\norm{ h(t,Y(t))}\wedge 1\bigr) d\mu(t)\leq
      2\epsilon.$$
This shows that $t\rightarrow h(t,Y(t))$ is
in $\ER(\R,\ellp{0}(\esprob,\prob,\espY),\mu)$.

For $p>0$, let $\delta>0$. From  the uniform integrability of
$\bigl(\norm{h(t,Y(t))}^p\bigr)_{t\in\R}$ 
(thanks to \eqref{eq:growth} and the uniform integrability of
$(\norm{Y(t)}^p)_{t\in\R}$), 
we can choose $\epsilon$ small
enough such that, for any measurable $A\subset\esprob$ such that $\prob(A)<\epsilon$,
$$\sup_{t\in\R}\expect\CCO{\un{A}\norm{h(t,Y(t)}^p}<\delta.$$
Note also that, for $p<1$, the mapping $U\mapsto\CCO{\expect\norm{U}^p}^{1/p}$, $\ellp{p}\rightarrow\R$, does
not satisfy the triangular inequality. However, the mapping $(U_1,U_2)\mapsto{\expect\norm{U_1-U_2}^p}$ is a
distance on $\ellp{p}$. We deduce that, for all $U_1,U_2,U_3\in\ellp{p}$,
$$\CCO{\expect\norm{U_1+U_2+U_3}^p}^{1/p}
\leq 3^{1/p-1}\left\lgroup\CCO{\expect\norm{U_1}^p}^{1/p} +\CCO{\expect\norm{U_2}^p}^{1/p}
+\CCO{\expect\norm{U_3}^p}^{1/p} \right\rgroup.
$$
To cover simultaneously the cases $p<1$ and $p\geq 1$,
we set $\kappa=\max(1,3^{1/p-1})$.
Using the same method as in the case when $p=0$, we get
 \begin{multline*}
\bigl(\expect\norm{ h(t,Y(t)) }^p)^{1/p}\\
\begin{aligned}
\leq&\kappa\biggl(\expect\bigl(\min_{1\leq i\leq m}
      \un{\Omega_{\epsilon,t}}  \norm{ h(t,Y(t))-h(t,y_i) }^p\bigr)\biggr)^{1/p}
     +\kappa\max_{1\leq i\leq m}\norm{h(t,y_i)}\\
 &+\kappa\biggl(\expect\bigl(\un{\Omega_{\epsilon,t}^c}\norm{h(t,Y(t))}^p\bigr)\biggr)^{1/p}\\
\leq&\kappa\Bigl(\max_{1\leq i\leq m}\norm{h(t,y_i)} +\epsilon+\delta\Bigr).
\end{aligned}
\end{multline*}
We conclude,  using  the ergodicity of $h(t,y_i)$ for all  $i\in \{1,\dots, m\}$,
that $t\rightarrow h(t,Y(t))$ is
in $\ER(\R,\ellp{p}(\esprob,\prob,\espY),\mu)$.

%%%%%
Secondly, we show that $F(.):=f(.,X(.))-f(.,Y(.))$ is in
$\ER(\R,\ellp{0}(\esprob,\prob,\espY),\mu)$. Let
$$\Phi :\,\left\{\begin{array}{lcl}
\espX&\rightarrow&\Cont(\R,\espY)\\
x&\mapsto&f(.,x).
\end{array}\right.$$

Let us endow $\Cont(\R,\espY)$ with the distance  $\dist_\infty(\varphi,\psi)
=\sup_{t\in\R}\norm{\varphi(t)-\psi(t)}$.
By the uniform
continuity assumption, $\Phi$ is continuous.
By \cite[Lemma 5.6]{blot-cieutat_ezzinbi2012}, for each
$\epsilon>0$, there exists $\eta>0$ such that,
for all $x,y\in\espX$,
\begin{equation}\label{eq:Schwartz}
 \Bigl\lgroup
   x \in K_{\epsilon} \text{ and } \dist(x,y)\leq \eta \Bigr\rgroup\
 \Rightarrow
\dist_\infty\Bigl(\Phi(x),\Phi(y)\Bigr)\leq\epsilon.
\end{equation}
For each $t\in\R$, let $\Omega_{\epsilon,t}$ be the subset of $\Omega$
on which $Y(t)\in K_{\epsilon}$. 
Since
$Z(t)=X(t)-Y(t)$,
we obtain, the following
inequalities, with the help of \eqref{eq:Schwartz} and
Chebyshev's inequality:
\begin{multline*}
\frac{\mu\accol{t\in[-r,r]\tq \expect
  \bigl(\norm{F(t)}\wedge 1\bigr)> 3\epsilon}}{\mu([-r,r])}\\
\begin{aligned}
\leq& \frac{\mu\accol{t\in[-r,r]\tq \expect \bigl(
      \un{\Omega_{\epsilon,t}}\un{\{\norm{Z(t)}>\eta\}}
             (\norm{F(t)}\wedge 1) \bigr)>
      \epsilon}}{\mu([-r,r])}\\
 &+\frac{\mu\accol{t\in[-r,r]\tq \expect \bigl(
      \un{\Omega_{\epsilon,t}}\un{\{\norm{Z(t)} \leq \eta\}}
             (\norm{F(t)}\wedge 1) \bigr)>
      \epsilon}}{\mu([-r,r])}\\
 &+\frac{\mu\accol{t\in[-r,r]\tq \expect \bigl(
      \un{\Omega^c_{\epsilon,t}} (\norm{F(t)}\wedge 1) \bigr)>
      \epsilon}}{\mu([-r,r])}\\
=& \frac{\mu\accol{t\in[-r,r]\tq \expect \bigl(
      \un{\Omega_{\epsilon,t}}\un{\{\norm{Z(t)} >\eta\}}
             (\norm{F(t)}\wedge 1) \bigr)>
      \epsilon}}{\mu([-r,r])}\\
 &+\frac{\mu\accol{t\in[-r,r]\tq \expect \bigl(
      \un{\Omega^c_{\epsilon,t}} (\norm{F(t)}\wedge 1) \bigr)>
      \epsilon}}{\mu([-r,r])}\\
\leq & \frac{\mu\accol{t\in[-r,r]\tq
         \prob\{  \norm{Z(t)}> \eta \}>\epsilon}}{\mu([-r,r])}\\
&+\frac{\mu\accol{t\in[-r,r]\tq \prob
    (\Omega^c_{\epsilon,t})>\epsilon}}{\mu([-r,r])}\\
\leq & \frac{\mu\accol{t\in[-r,r]\tq \frac{1}{ \eta}\expect \bigl(
       {\norm{Z(t)}} \bigr)> \epsilon}}{\mu([-r,r])}.
\end{aligned}
\end{multline*}
Since $Z$ is in $\ER(\R,\ellp{0}(\esprob,\prob,\espX),\mu)$, we have, for the above $\epsilon$,
$$\lim_{r\rightarrow \infty}\frac{\mu\accol{t\in[-r,r]\tq \expect \bigl(
       {\norm{Z(t)} \bigr)> \epsilon \eta}}}{\mu([-r,r])}=0,$$
which implies,  using \cite[Theorem 2.14]{blot-cieutat_ezzinbi2012}
(see Remark \ref{rem:ER_topological}),
$$\limsup_{r\rightarrow\infty}\frac{1}{\mu([-r,r])}\int_{[-r,r]}
      \expect\bigl(\norm{f(t,X(t))-f(t,Y(t))}\wedge 1\bigr)d\mu(t)=0.$$
Therefore $t\rightarrow f(t, X(t))-f(t,Y(t))$
is in $\ER(\R,\ellp{0}(\esprob,\prob,\espY),\mu)$.

%%%%%%%%%%%%%%%%%%%%%%%%%%%%%%%%%%%%%%%%%%%%%%%%%%%%%%%
Assume now that $p>0$. 
If $f$ is Lipschitz with respect to the second variable, uniformly
with respect to the first one, then $\norm{F(t)}\leq K\norm{Z(t)}$ for
some constant $K$, thus, trivially, 
$F\in \ER(\R,\ellp{p}(\esprob,\prob,\espY),\mu)$. 
% %%%
% If the family $(\norm{Z(t)}^p)_{t\in\R}$ is uniformly integrable, then 
% $(\norm{F(t)}^p)_{t\in\R}$ is also uniformly integrable. 
% Then, we can use the same reasoning as for $p=0$, replacing
% $\norm{F(t)}\wedge 1$ by $\norm{F(t)}^p$, since 
% $\expect \bigl(
%       \un{\Omega^c_{\epsilon,t}} (\norm{F(t)}^p) \bigr)$ 
% can be made arbitrarily small for $\epsilon$ sufficiently small. 
% We obtain
% \begin{eqnarray*}
%   \frac{\mu\accol{t\in[-r,r]\tq \expect
%   \bigl(\norm{F(t)}^p\bigr)> 3\epsilon}}{\mu([-r,r])} \leq
%   \frac{\mu\accol{t\in[-r,r]\tq \frac{1}{ \eta^p}\expect \bigl(
%        {\norm{Z(t)}^p} \bigr)> \epsilon}}{\mu([-r,r])}.
% \end{eqnarray*}
% Since $Z$ is in $\ER(\R,\ellp{p}(\esprob,\prob,\espX),\mu)$,
% we have, for the above $\epsilon$,
% $$\lim_{r\rightarrow \infty}\frac{\mu\accol{t\in[-r,r]\tq \expect \bigl(
%        {\norm{Z(t)}^p \bigr)> \epsilon \eta^p}}}{\mu([-r,r])}=0.$$
% We conclude, using \cite[Theorem 2.14]{blot-cieutat_ezzinbi2012}
% (see Remark \ref{rem:ER_topological}),
%  and the boundedness of $F(t)$ in
% $\ellp{p}(\esprob,\prob,\espY)$, that $t\rightarrow f(t, X(t))-f(t,Y(t))$ is in
% $\ER(\R,\ellp{p}(\esprob,\prob,\espY),\mu)$. 

\finproof
%%%%%%%%%

%%%%%%%%%%%%%%%%%%%%%%%%%%%%%%%%%%%%%%%%%%%%%%%%%%%%%%%%%%%%%%%%
\subsection{Pseudo-almost automorphy in $p$-mean
{vs} in $p$-dis\-tri\-bu\-tion}
Let $X=(X_t)_{t\in\R}$ be a continuous stochastic process with
values in $\espX$, defined on a probability space $(\esprob, \tribu,\prob)$. Let $\mu$ be a Borel measure on
$\R$ satisfying \eqref{eq:mu}. Clearly, we have for all $p\geq 0$,
$$\bigl\lgroup X\in  \AlA(\R,\ellp{p}(\esprob,\prob,\espX))\bigr\rgroup
\Rightarrow \bigl\lgroup X\in\AADf^p(\R,\espX)\bigr\rgroup.$$
Using Theorem \ref{theo:AAFD}, we can get more:
if $X$ satisfies \eqref{eq:Condition(2)-BMRF}, we deduce, for every
$p\geq 0$,
\begin{align*}
  \bigl\lgroup X\in \AlA(\R,\ellp{p}(\esprob,\prob,\espX)) \bigr\rgroup
  &\Rightarrow
  \bigl\lgroup X \in \AAD^p(\R,\espX)\bigr\rgroup,\\
  \bigl\lgroup X\in\PAA(\R,\ellp{p}(\esprob,\prob,\espX),\mu) \bigr\rgroup
  &\Rightarrow
  \bigl\lgroup X\in\PAAD^p(\R,\espX,\mu)\bigr\rgroup.
\end{align*}

The converse implications are false. Indeed,
Example \ref{exple:OU} shows that a process which is
almost automorphic in distribution is not necessarily
almost automorphic in probability or in $p$-mean, see also
\cite[Counterexample~2.16]{bedouhene-mellah-prf2012}\footnote{
{Let us here point out an infortunate error in
  \cite{MRF13}:
it is mistakenly said at the end of
  Section 1 of \cite{MRF13}
that almost periodicity in square mean implies almost periodicity in
distribution, and that the converse is true under a tightness
condition. The first claim is true under
the tightness condition \eqref{eq:Condition(2)-BMRF},
whereas Example \ref{exple:OU}, which is also Example 2.1 of
\cite{MRF13}, disproves the second claim.}
}.

The same counterexample also shows that a process which is
$\mu$-pseudo almost automorphic in $p$-distribution
is not necessarily
$\mu$-pseudo almost automorphic in probability or in
$p$-mean.

%%%%%%%%%%%%%%%%%%%%%%%%%%%%%%%%%%%%%% SDEs
\section{Pseudo almost automorphic solutions
to stochastic differential equations}\label{sec:SDE}

In the sequel, if $\espX$ and $\espY$ are metric spaces, we denote
$\CUB(\espX,\espY)$ the space of bounded uniformly continuous
functions from $\espX$ to $\espY$.

We are given two separable Hilbert spaces $\h_1$ and $\h_2$,
and
we consider the semilinear stochastic differential equation,
\begin{equation}\label{eq:SDE}
 dX_t = AX(t)\,dt + f(t, X(t))\,dt + g(t, X(t))\,dW(t), \ t\in\R
\end{equation}
where $A: \Dom(A)\subset\h_2\rightarrow \h_2$ is a densely defined closed
(possibly unbounded) linear operator, and
$f : \R\times\h_2 \rightarrow\h_2$,
and $g : \R\times\h_2\rightarrow L(\h_1, \h_2)$
are continuous functions.
In this section, we assume that:

\newcounter{carbone}
\newcounter{carbon2}
\begin{enumerate}[{\rm (i)}]

 \item \label{cond:wienerpro}
 $W(t)$ is an $\h_1$-valued Wiener process with nuclear covariance operator
 $Q$ (we denote by $\trace Q$ the trace of $Q$),
 defined on a stochastic basis
 $(\esprob,\tribu,(\tribu_t)_{t\in \R},\prob)$.

\setcounter{carbon2}{\value{enumi}}
 \item \label{cond:semigrp} $A : \Dom(A)\rightarrow \h_2$ is the
infinitesimal generator of
a $C_0$-semigroup $(S(t))_{t\geq0}$ such that there exists a constant
$\delta>0$ with
$$\|S(t)\|_{L(\h_2)}\leq e^{-\delta t}, t\geq0.$$
  \item \label{cond:croissance}
There exists a constant $K$ such that the mappings
$f :\R\times\h_2\rightarrow \h_2$ and
$g :\R\times\h_2\rightarrow L(\h_1, \h_2)$
satisfy
$$\|f(t,x)\|_{\h_{2}} + \|g(t,x)\|_{L(\h_1, \h_2)} \leq K (1+\|x\|_{\h_2}) .$$

\item \label{cond:lipschitz}The functions $f$ and $g$ are Lipschitz,
  more precisely there exists a constant $K$ such that
$$\|f(t,x) - f(t,y)\|_{\h_2} + \|g(t,x) - g(t,y)\|_{L(\h_1, \h_2)}\leq
K\|x - y\|_{\h_2}$$ for all $t\in \R$
and $x, y \in \h_2$.
\setcounter{carbone}{\value{enumi}}
\item \label{cond:PAA}
$f \in \PAAUb(\R\times \h_2, \h_2,\mu) $ and
$g \in \PAAUb(\R\times \h_2, L(\h_1, \h_2),\mu)$
for some given
Borel measure $\mu$ on $\R$ which satisfies \eqref{eq:mu}
and Condition (\textbf{H}).

\end{enumerate}
By \cite[Theorem 3.5]{blot-cieutat_ezzinbi2012},
Condition (\textbf{H}) implies that $\ER(\R,\espX,\mu)$ and
$\PAA(\R,\espX,\mu)$ are translation invariant.

In order to study the
weighted pseudo almost automorphy property of
solutions of SDEs, we need a result on  almost automorphy.
%%%%%%%%%%%%%%%%%%%%%%%%%%%%%%
\begin{theorem}\label{theo:AA}
{\bf(Almost automorphic solution of an equation with almost
  automorphic coefficients)}
 Let the assumptions (\ref{cond:wienerpro}) - (\ref{cond:lipschitz}) be
 fulfilled, and assume furthermore the following condition, which is
 stronger than (\ref{cond:PAA}) :
\begin{itemize}
\item[(\ref{cond:PAA}')] $f \in \AlAUb(\R\times \h_2, \h_2) $ and
$g \in \AlAUb(\R\times \h_2, L(\h_1, \h_2))$.
\end{itemize}
Assume further
that
 $\theta :=\dfrac{K^2}{\delta}\CCO{\dfrac{1}{2\delta} + \trace Q}<1$.
Then there exists a unique mild solution $X$ to (\ref{eq:SDE})
in the space $\CUB\bigl(\R,\ellp{2}(\prob, \h_2)\bigr)$ of bounded
uniformly continuous mappings from $\R$ to $\ellp{2}(\prob, \h_2)$.
Furthermore,
$X$ has a.e.~continuous
trajectories,
and $X(t)$ satisfies, for each $t\in \R$:
\begin{equation}\label{eq:mildsol}
X(t) = \int^{t}_{-\infty}S(t-s)f\bigl(s, X(s)\bigr)ds +
\int^{t}_{-\infty}S(t-s)g\bigl(s, X(s)\bigr)dW(s).
\end{equation}
If furthermore
$\theta' :=\dfrac{4K^2}{\delta}\CCO{\dfrac{1}{\delta} + \trace Q}<1$,
then $X$ is almost
automorphic in 2-distribution.
\end{theorem}
%%%%%%%%%%%%

The proof of this theorem is very similar to that
of \cite[Theorem 3.1]{KMRF12averaging},
which is the analogous result
for SDEs with almost periodic coefficients.
%(a preliminary version of this proof can be found in \cite{omarPHD}).
Only the almost automorphy part needs to be adapted.
Such an adaptation is provided in \cite{liu-sun2014},
for SDEs driven by L\'evy processes, but only for
one-dimensional almost automorphy.
We give the proof of this part for the convenience of the
reader.

Let us first recall the following result, which is given in a more general
form in \cite{DaPrato-Tudor95}:
%%%% DA PRATO-TUDOR Prop. 3.1.(c) %%%%
\begin{proposition}\label{prop:daprato-tudor3.1}%
(\cite[Proposition 3.1-(c)]{DaPrato-Tudor95})
Let $\tau\in\R$.
Let $(\xi_n)_{0\leq n\leq\infty}$
be a sequence of square integrable $\h_2$-valued random variables.
Let $(f_n)_{0\leq n\leq\infty}$
and $(g_n)_{0\leq n\leq\infty}$ be sequences of mappings from $\R\times\h_2$
to $\h_2$ and $L(\h_1, \h_2)$ respectively, satisfying
\eqref{cond:croissance} and \eqref{cond:lipschitz} (replacing $f$ and
$g$ by $f_n$ and $g_n$ respectively, and the constant $K$ being
independent of $n$).
For each $n$, let $X_n$ denote the solution to
\begin{multline*}
X_n(t)=S(t-\tau)\xi_n\\
+\int^{t}_{\tau}S(t-s)f_n\bigl(s, X_n(s)\bigr)ds +
\int^{t}_{\tau}S(t-s)g_n\bigl(s, X_n(s)\bigr)dW(s).
\end{multline*}
Assume that, for every $(t,x)\in \R\times\h_2$,
\begin{gather*}
\lim_{n\rightarrow\infty}f_n(t,x)=f_{\infty}(t,x),\
\lim_{n\rightarrow\infty}g_n(t,x)=g_{\infty}(t,x),\\
\lim_{n\rightarrow\infty}\dist_{\bl}(\law{\xi_n,W},\law{\xi_\infty,W})=0,
\end{gather*}
(the last equality takes place in $\laws{\h_2\times \Cont(\R,\h_1)}$). Then we have in $\Cont([\tau,T]; \h_2)$,
for any $T>\tau$,
$$\lim_{n\rightarrow\infty}\dist_{\bl}(\law{X_n},\law{X_\infty})=0.$$
\end{proposition}

%%% Gronwall
We need also a variant of Gronwall's lemma.
%%%%%%%%%%%%%%%%%%%%%%%%%%%%%%%
\begin{lemma}\label{lem:gronwall}
(\cite[Lemma 3.3]{KMRF12averaging})
Let $g :\R\to\R$ be a continuous function
such that,
for every $t\in\R$,
\begin{equation}\label{hypothesegronwall}
0\leq g(t)\leq \alpha(t)
+\beta_1\int_{-\infty}^t e^{-\delta_1(t-s)}g(s)\,ds
+\dots
+\beta_n\int_{-\infty}^t e^{-\delta_n(t-s)}g(s)\,ds,
\end{equation}
for some locally integrable function $\alpha :\,\R\rightarrow\R$,
and for some constants $\beta_1,\dots,\beta_n\geq 0$,
and some constants $\delta_1,\dots,\delta_n>\beta$,
where $\beta:=\sum_{i=1}^n\beta_i$.
We assume that the integrals in the right hand
side of \eqref{hypothesegronwall} are convergent.
Let $\delta=\min_{1\leq i\leq n}\delta_i$.
Then, for every
$\gamma\in ]0, \delta-\beta]$
such that $\int_{-\infty}^0 e^{\gamma s}\alpha(s)\,ds$ converges, we
have, for every $t\in\R$,
\begin{equation*}%\label{eq:gronwallinfinite}
g(t)\leq \alpha(t)+\beta\int_{-\infty}^t e^{-\gamma(t-s)}\alpha(s)\,ds.
\end{equation*}
In particular, if $\alpha$ is constant, we have
\begin{equation*}%\label{eq:alphaconstant}
g(t)\leq \alpha\, \frac{\delta}{\delta-\beta}.
\end{equation*}
\end{lemma}
%%%%%%%%%

%%%%%%%%%%%%%%%%%%%%%%%%%%%%%%%
\proofof{Theorem \ref{theo:AA}}
The proof of the existence and uniqueness
of a mild solution to (\ref{eq:SDE}) in
$\CUB\bigl(\R,\ellp{2}(\prob, \h_2)\bigr)$
is the same as that of Theorem 3.1 in
\cite{KMRF12averaging} or Theorem 3.3.1 in \cite{omarPHD}.

For the almost automorphy part, let $(\gamma'_n)$ be a sequence in $\R$. Since $f$ and $g$ are almost
automorphic, there exists a subsequence $(\gamma_n)$ and functions $\limaa{f} :\,\R\times \h_2\rightarrow \h_2$
and $\limaa{g} :\,\R\times \h_2\rightarrow L(\h_1,\h_2)$ such that
\begin{align}\label{eq:APF}
 \lim_{n\rightarrow\infty}f(t+\gamma_n, x)&= \limaa{f}(t,x), \quad
 \lim_{n\rightarrow\infty}\limaa{f}(t-\gamma_n, x)= {f}(t,x)\\
\label{eq:APG}
 \lim_{n\rightarrow\infty}g(t+\gamma_n,x)&= \limaa{g}(t,x),\quad
\lim_{n\rightarrow\infty}\limaa{g}(t-\gamma_n, x)= {g}(t,x).
\end{align}
These limits are taken uniformly with respect to $x$ in
bounded subsets of $\h_2$.

For each fixed integer $n$, we consider
$$X_n(t) = \int^t_{-\infty}S(t-s)f(s+\gamma_n, X_n(s))\,ds
     + \int^t_{-\infty}S(t-s)g(s+\gamma_n, X_n(s))\,dW(s)$$
the mild solution to
\begin{equation*}
 dX_n(t) = AX_n(t)dt + f(t+\gamma_n, X_n(t))\,dt + g(t+\gamma_n, X_n(t))\,dW(t)
\end{equation*}
and
$$\limaa{X}(t)
= \int^t_{-\infty}S(t-s)\limaa{f}(s, \limaa{X}(s))\,ds
+ \int^t_{-\infty}S(t-s)\limaa{g}(s, \limaa{X}(s))\,dW(s)$$
the mild solution to
\begin{equation*}
 d\limaa{X}(t) = A(t)\limaa{X}(t)dt + \limaa{f}(t, \limaa{X}(t))\,dt
                + \limaa{g}(t, \limaa{X}(t))\,dW(t).
\end{equation*}
Make the change of variable $\sigma + \gamma_n = s$, the process
\begin{multline*}
X(t+ \gamma_n)
= \int^{t+\gamma_n}_{-\infty}S(t+\gamma_n-s)f(s, X(s))\,ds\\
     + \int^{t+\gamma_n}_{-\infty}S(t+\gamma_n-s)g(s, X(s))\,dW(s)
\end{multline*}
satisfies
\begin{multline*}
X(t+ \gamma_n)= \int^t_{-\infty}S(t-s)f(s+\gamma_n, X(s+\gamma_n))ds \\
     + \int^t_{-\infty}S(t-s)g(s+\gamma_n,X(s+\gamma_n))d\tilde{W}_n(s),
\end{multline*}
where  $\tilde{W}_n(s) = W(s+\gamma_n) - W(\gamma_n)$
is a Brownian motion with the same distribution as $W(s)$.
Thus the process $X(.+\gamma_n)$ has the same distribution as $X_n$.

Let us show that $X_n(t)$ converges in quadratic
mean to $\limaa{X}(t)$ for each fixed
$t\in \R$.
We have
\begin{multline*}
 \expect\lVert X_n(t) - \limaa{X}(t)\rVert^2 \\
\begin{aligned}
=&
\expect\biggl\lVert \int^t_{-\infty}S(t-s)
  \Bigl( f(s+\gamma_n, X_n(s)) - \limaa{f}(s, \limaa{X}(s))\Bigr)\,ds\\
&+\int^t_{-\infty}S(t-s)\Bigl( g(s+\gamma_n, X^n(s))
      - \limaa{g}(s, \limaa{X}(s))\Bigr)\,dW(s)
       \biggr\rVert^2\\
 \leq&2\expect\biggl\lVert \int^t_{-\infty}S(t-s)
    \Bigl( f(s+\gamma_n, X_n(s))
                      - \limaa{f}(s, \limaa{X}(s))\Bigr)\,ds\biggr\rVert^2\\
&+2\expect\biggl\lVert
\int^t_{-\infty}S(t-s)
   \Bigl(g(s+\gamma_n, X_n(s)) - \limaa{g}(s,\limaa{X}(s))\Bigr)\,dW(s)
       \biggl\rVert^2\\
 \leq&4\expect\biggl\lVert \int^t_{-\infty}S(t-s)
      \Bigl( f(s+\gamma_n, X^n(s))
               -f(s+\gamma_n, \limaa{X}(s))\Bigr)\,ds
       \biggr\rVert^2\\
&+4\expect\biggl\lVert\int^t_{-\infty}S(t-s)
      \Bigl( f(s+\gamma_n, \limaa{X}(s))- \limaa{f}(s, \limaa{X}(s))\Bigr)\,ds
      \biggr\rVert^2\\
&+4\expect\biggl\lVert
\int^t_{-\infty}S(t-s)
      \Bigl( g(s+\gamma_n, X^n(s)) -g(s+\gamma_n, \limaa{X}(s))\Bigr)\,dW(s)
      \biggr\rVert^2\\
&+4\expect\biggl\lVert
\int^t_{-\infty}S(t-s)
      \Bigl( g(s+\gamma_n, \limaa{X}(s))-\limaa{g}(s, \limaa{X}(s))\Bigr)\,dW(s)
        \biggr\rVert^2\\
\leq& I_1 + I_2 + I_3+ I_4.
\end{aligned}
\end{multline*}
Now, using (\ref{cond:semigrp}), (\ref{cond:lipschitz})
and the Cauchy-Schwartz inequality, we obtain
\begin{align*}
I_1&= 4\expect\biggl\lVert \int^t_{-\infty}S(t-s)
      \Bigl( f(s+\gamma_n, X_n(s)) -f(s+\gamma_n, \limaa{X}(s))\Bigr)\,ds
              \biggr\rVert^2\\
&\leq4\expect\biggl(\int^t_{-\infty}\lVert
  S(t-s)\rVert \lVert f(s+\gamma_n, X_n(s)) -f(s+\gamma_n, \limaa{X}(s))
  \rVert \,ds\biggr)^2\\
&\leq4\expect\biggl(\int^t_{-\infty}e^{-\delta(t-s)}
  \lVert f(s+\gamma_n, X_n(s)) -f(s+\gamma_n, \limaa{X}(s))\rVert ds\biggr)^2\\
&\leq4\biggl(\int^t_{-\infty}e^{-\delta(t-s)}\,ds\biggr)
    \biggl(\int^t_{-\infty}e^{-\delta(t-s)}\expect\lVert f(s+\gamma_n, X_n(s))
               -f(s+\gamma_n, \limaa{X}(s))\rVert^2 ds\biggr)\\
&\leq \frac{4K^2}{\delta}\int^t_{-\infty}e^{-\delta(t-s)}
                                \expect\lVert X_n(s) - \limaa{X}(s)\rVert^2ds.
\end{align*}
Then we have
\begin{align*}
 I_2 & = 4\expect\biggl\lVert\int^t_{-\infty}S(t-s)[f(s+\gamma_n, \limaa{X}(s))
          - \limaa{f}(s, \limaa{X}(s))]ds\biggr\rVert^2\\
&\leq 4\expect\biggl(\int^t_{-\infty}e^{-\delta(t-s)}
     \lVert f(s+\gamma_n, \limaa{X}(s)) -\limaa{f}(s, \limaa{X}(s))\rVert ds\biggr)^2\\
&\leq 4\expect\biggl(\int^t_{-\infty}e^{-\delta(t-s)}ds\biggr)
\biggl(\int^t_{-\infty}e^{-\delta(t-s)}
         \lVert f(s+\gamma_n, \limaa{X}(s)) -\limaa{f}(s, \limaa{X}(s))\rVert^2 ds\biggl)\\
&\leq4 \biggl(\int^t_{-\infty}e^{-\delta(t-s)}ds\biggr)^2
        \sup_{s}\expect\lVert f(s+\gamma_n, \limaa{X}(s)) -\limaa{f}(s, \limaa{X}(s))\rVert^2\\
&\leq \frac{4}{\delta^2}\sup_{s}
             \expect\lVert f(s+\gamma_n, \limaa{X}(s)) -\limaa{f}(s, \limaa{X}(s))\rVert^2,
\end{align*}
which converges to $0$ as $n\rightarrow\infty$
because $\sup_{t\in\R} \expect\lVert \limaa{X}(t)\rVert^2 < \infty$
which implies that $(\limaa{X}(t))_{t}$ is tight relatively to bounded sets.

Applying It\^o's isometry, we get
\begin{align*}
 I_3& = 4\expect\biggl\lVert\int^t_{-\infty}S(t-s)
   \Bigl(g(s+\gamma_n, X_n(s)) -g(s+\gamma_n, \limaa{X}(s))\Bigr)\,
                 dW(s)\biggr\rVert^2\\
& \leq 4\trace Q\expect\int^t_{-\infty}\lVert S(t-s)\rVert^2 \,
      \lVert g(s+\gamma_n, X_n(s)) -g(s+\gamma_n, \limaa{X}(s))\rVert^2 ds\\
&\leq 4\trace Q\int^t_{-\infty}e^{-2\delta(t-s)} \expect\lVert
        g(s+\gamma_n, X_n(s)) -g(s+\gamma_n, \limaa{X}(s))\rVert^2 ds\\
&\leq 4K^2\trace Q\int^t_{-\infty}e^{-2\delta(t-s)}
      \expect\lVert X_n(s) - \limaa{X}(s)\rVert^2ds,
\end{align*}
and
\begin{align*}
 I_4& = 4\expect\biggl\lVert\int^t_{-\infty}S(t-s)
     \Bigl(g(s+\gamma_n, \limaa{X}(s))-\limaa{g}(s, \limaa{X}(s))\Bigr)
        \, dW(s)\biggr\rVert^2\\
& \leq 4\trace Q\expect\biggl(\int^t_{-\infty}\lVert S(t-s)\rVert^2
      \lVert g(s+\gamma_n, \limaa{X}(s))-\limaa{g}(s, \limaa{X}(s))\rVert^2 ds\biggr)\\
&\leq 4\trace Q\biggl(\int^t_{-\infty}e^{-2\delta(t-s)}ds\biggr)
      \sup_{s\in\R} \expect\lVert g(s+\gamma_n, \limaa{X}(s))
                  -\limaa{g}(s, \limaa{X}(s))\rVert^2\\
& \leq \frac{2\trace Q}{\delta}\,
     \sup_{s\in\R} \expect\lVert g(s+\gamma_n, \limaa{X}(s))
          -\limaa{g}(s, \limaa{X}(s))\rVert^2.
\end{align*}
For the same reason as for $I_2$,
the right hand term goes to $0$ as
$n\rightarrow\infty$.

We thus have
\begin{multline*}
\expect\lVert X_n(t) - \limaa{X}(t)\rVert^2
\leq \alpha_n
  +\frac{4K^2}{\delta}\int^t_{-\infty}e^{-\delta(t-s)} \expect\lVert
  X_n(s)-\limaa{X}(s)\rVert^2\,ds\\
  +4K^2\trace{Q}\int^t_{-\infty}e^{-2\delta(t-s)} \expect\lVert
  X_n(s)-\limaa{X}(s)\rVert^2\,ds
\end{multline*}
for a sequence $(\alpha_n)$ such that
$\lim_{n\rightarrow\infty}\alpha_n=0.$ Furthermore,
$\beta:=\frac{4K^2}{\delta}+4K^2\trace{Q}<\delta$. We conclude by
 Lemma \ref{lem:gronwall} that
$$\lim_{n\rightarrow\infty}  \expect\lVert X_n(t) - \limaa{X}(t)\rVert^2 = 0,$$
hence $X_n(t)$ converges in distribution to $\limaa{X}(t)$.
But, since the distribution of $X_n(t)$ is the same as that of
$X(t+\gamma_n)$,
we deduce that $X(t+\gamma_n)$ converges in distribution to
$\limaa{X}(t)$.
By analogy and using (\ref{eq:APF}), (\ref{eq:APG}) we can easily prove that
$\limaa{X}(t-\gamma_n)$ converges in distribution to
${X}(t)$.

Note that the sequence $(\norm{X_n(t)}^2)$ is uniformly integrable, thus
$(\norm{X(t+\gamma_n)}^2)$ is uniformly integrable too.
As $(\gamma'_n)$ is arbitrary, this implies that the family
$(\norm{X(t)}^2)_{t\in\R}$ is uniformly integrable, because, if not,
there would exist a sequence $(\gamma'_n)$ and $t\in\R$ such that
no subsequence of $(\norm{X(t+\gamma'_n)}^2)$ is uniformly integrable.

We have thus proved that $X$
has almost automorphic one-dimensional 2-distributions.
To prove that $X$ is almost automorphic in 2-distribution,
we apply
Proposition \ref{prop:daprato-tudor3.1}: for fixed $\tau\in\R$, let
$\xi_n=X(\tau+\gamma_n)$, $f_n(t,x)=f(t+\gamma_n,x)$,
$g_n(t,x)=g(t+\gamma_n,x)$.
By the foregoing, $(\xi_n)$ converges in distribution to some variable
$Y(\tau)$.
We deduce that $(\xi_n)$ is tight,
and thus $(\xi_n,W)$ is tight also.
We can thus choose a subsequence (still noted $(\gamma_n)$ for
simplicity)
such that $(\xi_n,W)$ converges in
distribution to $(Y(\tau),W)$.
Then, by Proposition \ref{prop:daprato-tudor3.1}, for every $T\geq \tau$,
$X(.+\gamma_n)$ converges in distribution on $\Cont([\tau,T]; \h_2)$ to
the (unique in distribution) solution to
\begin{equation*}
Y(t)=S(t-\tau)Y(\tau)+\int^{t}_{\tau}S(t-s)f\bigl(s, Y(s)\bigr)\,ds
+ \int^{t}_{\tau}S(t-s)g\bigl(s, Y(s)\bigr)\,dW(s).
\end{equation*}
Note that $Y$ does not depend on the chosen interval $[\tau,T]$, thus
the convergence takes place on $\Cont(\R; \h_2)$.
Similarly, $Y_n:=Y(.-\gamma_n)$ converges in distribution on
$\Cont(\R; \h_2)$ to $X$.
Thus $X$ is almost automorphic in 
2-distribution.
\finproof
%%%%%%%%%

We are now ready to prove our main result.
%%%%%%%%%%%%%%%%%%%%%%%%%%%%%%%%
\begin{theorem}\label{theo:main}
{\bf(Weighted pseudo almost automorphic solution of an equation with
  weighted pseudo almost automorphic coefficients)}
 Let the assumptions (\ref{cond:wienerpro}) - (\ref{cond:PAA}) be
 fulfilled.
Let $(f_1,g_1)$ and $(f_2,g_2)$ be respectively the decompositions of $f$ and $g$, namely,
\begin{gather*}
f=f_1+f_2, \quad g=g_1+g_2, \\
f_1\in\AlAUb(\R\times \h_2, \h_2),\
f_2\in \ERUb(\R\times \h_2, \h_2,\mu),\\
g_1\in\AlAUb(\R\times \h_2, L(\h_1, \h_2)),\
g_2\in \ERUb(\R\times \h_2, L(\h_1, \h_2),\mu).
\end{gather*}
Assume that $f_1$ and $g_1$ satisfy the same growth and Lipschitz conditions
(\ref{cond:croissance}) - (\ref{cond:PAA}) as $f$ and $g$
respectively, with same coefficient $K$.
Assume furthermore that
$$\theta' :=\dfrac{4K^2}{\delta}\CCO{\dfrac{1}{\delta} + \trace Q}<1.$$
%%%
Then there exists a unique mild solution $X$ to \eqref{eq:SDE}
in the space
$\CUB\bigl(\R,\ellp{2}(\prob, \h_2)\bigr)$ of bounded
uniformly continuous mappings from $\R$ to $\ellp{2}(\prob, \h_2)$,
$X$ has a.e.~continuous
trajectories,
and
$X$ satisfies \eqref{eq:mildsol}  for every $t\in \R$.
Furthermore, $X$ is $\mu$-pseudo almost automorphic in $2$-distribution.
More precisely, let
$Y\in\CUB\bigl(\R,\ellp{2}(\prob, \h_2)\bigr)$ be the unique
almost automorphic in distribution
mild solution to
\begin{equation}\label{eq:edsY}
dY(t)=
AY(t)\,dt + f_1(t, Y(t))\,dt + g_1(t, Y(t))\,dW(t), \ t\in\R.
\end{equation}
Then $X$ has the decomposition
$$X=Y+Z, \quad
Z\in\ER\bigl(\R,\ellp{2}(\prob,\h_2),\mu\bigr).
$$
\end{theorem}
%%%%%%%%%%%%%
The following technical lemma will be used several times.
%%%%%%%%%%%%%
\begin{lemma}\label{lem:convolut}
Let $\hm\in\ER(\R,\R,\mu)$. Then the function
$$t\mapsto\biggl(\int_{-\infty}^t e^{-2\delta(t-s)}\hm^2(s)\,ds\biggr)^{1/2}$$
is also in $\ER(\R,\R,\mu)$.
\end{lemma}
%%%%%%%%%%%
\proof
by Condition (\textbf{H}) and \cite[Theorem 3.9]{blot-cieutat_ezzinbi2012},
we have, for every $u\in\R$,
\begin{equation*}
\lim_{r\rightarrow+\infty}
 \frac{1}{\mu([-r,r])} \int_{[-r,r]}\abs{\hm(t-u)}\,d\mu(t)
=0.
\end{equation*}
We deduce,
by Lebesgue's dominated convergence theorem,
\begin{multline*}
 \frac{1}{\mu([-r,r])} \int_{[-r,r]} \biggl(\int_{-\infty}^t e^{-2\delta (t-s)}
            \hm^2(s)\, ds
        \biggr)^{1/2} d\mu(t)\\
\begin{aligned}
&\leq  \frac{1}{\big(\mu([-r,r])\big)^{1/2}}
   \biggl(\int_{[-r,r]} \int_{-\infty}^t e^{-2\delta (t-s)}
            \hm^2(s)\, ds
         \,d\mu(t)\biggr)^{1/2}\\
&=  \frac{1}{\big(\mu([-r,r])\big)^{1/2}}
    \biggl(\int_{[-r,r]} \int_{0}^{+\infty} e^{-2\delta u}
            \hm^2(t-u)\, du
         \,d\mu(t)\biggr)^{1/2}\\
&=  \frac{1}{\big(\mu([-r,r])\big)^{1/2}}
    \biggl( \int_{0}^{+\infty} e^{-2\delta u}
            \int_{[-r,r]}\hm^2(t-u)\,d\mu(t)\, du
         \biggr)^{1/2}\\
& \leq  \biggl(\int_{0}^{+\infty} e^{-2\delta u}
      %\frac{\norm{\hm_\epsilon}_\infty}{\mu([-r,r])}\,
      \norm{\hm}_{\infty}\,
       \frac{\int_{[-r,r]}{\abs{\hm(t-u)}}\,d\mu(t)}{\mu([-r,r])}\,du
      \biggr)^{1/2}\\
&\rightarrow 0 \text{ when }r\rightarrow +\infty.
\end{aligned}
\end{multline*}
\finproof
%%%%%%%%%

%%%%%%%%%%%%%%%%%%%%%%%%%%%%%%%%%%
\proofof{Theorem \ref{theo:main}}
The existence and the properties of $Y$ are guaranteed by Theorem
\ref{theo:AA}.

As in Theorem \ref{theo:AA},
the existence and uniqueness of the mild solution $X$ to \eqref{eq:SDE}
are proved as in \cite[Theorem 3.1]{KMRF12averaging},
 using the classical method of the fixed point theorem for the
contractive operator $L$ on
$\CUB\bigl(\R,\ellp{2}(\prob, \h_2)\bigr)$
defined by
$$LX(t) = \int^{t}_{-\infty}S(t-s)f\bigl(s, X(s)\bigr)ds +
\int^{t}_{-\infty}S(t-s)g\bigl(s, X(s)\bigr)dW(s).$$

The solution $X$ defined by \eqref{eq:mildsol} is thus the limit in $\CUB\bigl(\R,\ellp{2}(\prob, \h_2)\bigr)$
of a sequence $(X_n)$ with  arbitrary $X_0$ and, for every $n$, $X_{n+1}=L(X_n)$. To prove that $X$ is
$\mu$-pseudo almost automorphic in $2$-distribution we choose a special sequence. Set
$$
X_0=Y, \
X_{n+1}=L(X_n),\ Z_n=X_n-Y,\ n\in\N.
$$
Let us prove that each $Z_n$ is in
$\ER\bigl(\R,\ellp{2}(\prob,\h_2),\mu\bigr)$.
We use some arguments of the proof of
\cite[Theorem 5.7]{blot-cieutat_ezzinbi2012}.
We have, for every $n\in\N$ and every $t\in\R$,
\begin{align*}
Z_{n+1}(t)&=LX_n(t)-Y(t)\\
=&\int_{-\infty}^t S(t-s)\bigl( f(s,X_n(s))-f(s,Y(s))\bigr)\,ds\\
&+\int_{-\infty}^t S(t-s)\bigl( g(s,X_n(s))-g(s,Y(s))\bigr)\,dW(s)\\
&+\int_{-\infty}^t S(t-s)\bigl( f(s,Y(s))-f_1(s,Y(s))\bigr)\,ds\\
&+\int_{-\infty}^t S(t-s)\bigl( g(s,Y(s))-g_1(s,Y(s))\bigr)\,dW(s)\\
=&\int_{-\infty}^t S(t-s)\bigl( f(s,X_n(s))-f(s,Y(s))\bigr)\,ds\\
&+\int_{-\infty}^t S(t-s)\bigl( g(s,X_n(s))-g(s,Y(s))\bigr)\,dW(s)\\
&+\int_{-\infty}^t S(t-s)f_2(s,Y(s))\,ds +\int_{-\infty}^t S(t-s)g_2(s,Y(s))\,dW(s).
\end{align*}
Assume that $Z_n\in\ER\bigl(\R,\ellp{2}(\prob,\h_2),\mu\bigr)$.
By the Lipschitz condition \eqref{cond:lipschitz},
$$\bigl(\expect\norm{f(t,X_n(t))-f(t,Y(t))}^2\bigr)^{1/2}
\leq K\bigl(\expect\norm{Z_n(t)}^2\bigr)^{1/2}$$ thus the mapping
$$\fm : t\mapsto \bigl(\expect\norm{f(t,X_n(t))-f(t,Y(t))}^2\bigr)^{1/2}$$
is in $\ER(\R,\R,\mu)$. The same conclusion holds for
$$\gm : t\mapsto \bigl(\expect\norm{g(t,X_n(t))-g(t,Y(t))}^2\bigr)^{1/2}.$$
We get, using Lemma \ref{lem:convolut},
\begin{multline*}
\frac{1}{\mu([-r,r])}\int_{[-r,r]}
  \biggl(\expect \norm{
    {\int_{-\infty}^t S(t-s)\bigl( f(s,X_n(s))-f(s,Y(s))\bigr)\,ds}}^2
     \biggr)^{1/2} d\mu(t)\\
\begin{aligned}
&\leq  \frac{1}{\mu([-r,r])}\int_{[-r,r]}
           \CCO{\int_{-\infty}^t e^{-2\delta (t-s)}\fm^2(s) \,ds
        }^{1/2}d\mu(t)\\
&\rightarrow 0\text{ when }r\rightarrow +\infty,
\end{aligned}
\end{multline*}
and
\begin{multline*}
\frac{1}{\mu([-r,r])}\int_{[-r,r]}
  \biggl(\expect\norm{
    {\int_{-\infty}^t S(t-s)\bigl( g(s,X_n(s))-g(s,Y(s))\bigr)\,dW(s)}}^2
     \biggr)^{1/2} d\mu(t)\\
\begin{aligned}
&\leq  (\trace Q)^{1/2}\frac{1}{\mu([-r,r])}\int_{[-r,r]}
           \CCO{\int_{-\infty}^t e^{-2\delta (t-s)}\gm^2(s) \,ds
        }^{1/2}d\mu(t)\\
&\rightarrow 0\text{ when }r\rightarrow +\infty.
\end{aligned}
\end{multline*}

To prove that $Z_{n+1}$ is in
$\ER\bigl(\R,\ellp{2}(\prob,\h_2),\mu\bigr)$,
there only remains to show that the process
$\int_{-\infty}^t S(t-s)f_2(s,Y(s))\,ds
+\int_{-\infty}^t S(t-s)g_2(s,Y(s))\,dW(s)$
belongs to
$\ER\bigl(\R,\ellp{2}(\prob,\h_2),\mu\bigr)$.
As $Y$ is almost automorphic in distribution, the family
$(\transl{Y}(t))=({Y}(t+.))_{t\in\R}$ is uniformly tight in
$\Cont_k(\R,\h_2)$. In particular, for each $\epsilon>0$
there exists a compact subset $\mathcal{K}_\epsilon$ of
$\Cont_k(\R,\h_2)$ such that, for every $t\in\R$,
$$\prob\accol{\transl{Y}(t)\in\mathcal{K}_\epsilon}\geq 1-\epsilon.$$
%%%%%%%%%%%%%%%%%%%%%%%%%%%%%%%%%%%%%%%%%%%%%%%%%%%%%
By the Arzel\`a-Ascoli Theorem
(e.g.~\cite[Theorems 8.2.10 and 8.2.11]{Engelking89}),
this implies that,
for every $\epsilon>0$, and for every compact interval $I$ of $\R$,
there exists a compact subset $K_{\epsilon,I}$ of $\h_2$ such that,
for every $t\in\R$,
$$\prob\accol{(\forall s\in I)\ Y(t+s)\in K_{\epsilon,I}}\geq
1-\epsilon. $$
In particular, the family $(Y(t))_{t\in\R}$ is tight, i.e., 
denoting $K_{\epsilon}=K_{\epsilon,\{0\}}$, we have, for every
$t\in\R$, 
\begin{equation*}%\label{eq:Ytendu}
\prob\accol{Y(t)\in K_{\epsilon}}\geq 1-\epsilon.
\end{equation*}
By the uniform continuity property
of $f_2$ and $g_2$ on $K_{\epsilon}$, there exists
$\eta(\epsilon)>0$ such that, for all $x,y\in K_{\epsilon}$,
$$
\norm{x-y}\leq \eta(\epsilon)\Rightarrow
\sup_{t\in\R}\max\bigl(\norm{f_2(t,x)-f_2(t,y)},\norm{g_2(t,x)-g_2(t,y)}\bigr)
\leq \epsilon.$$
We can find a finite sequence
$y_1,\dots,y_{m}$ such that
$$K_{\epsilon}\subset
\bigcup_{i=1}^{m}B(y_i,\eta(\epsilon)).$$
By \cite[Remark 3.6]{KMRF12averaging}),
the condition $\theta'<1$ ensures that
$Y$  is bounded in $\ellp{p}(\prob, \h_2)$ for some $p>2$
(the same result holds for $X$, but we do not need it).
Note that $f_2=f-f_1$ and $g_2=g-g_1$ satisfy a condition similar to
\eqref{cond:croissance}, which implies that
$f_2(.,Y(.))$ and $g_2(.,Y(.))$ are bounded in $\ellp{p}(\prob,
\h_2)$ and $\ellp{p}\bigl(\prob, L(\h_1, \h_2)\bigr)$ respectively.
Let
$$\Mp=\sup_{t\in\R}\max\bigl( \expect\norm{f_2(.,Y(.))}^p,
                    \expect\norm{g_2(.,Y(.))}^p  \bigr)^{2/p}.$$
Let $q=p/(p-2)$.
Let $t\in\R$, and let 
$\Omega_{\epsilon,t}$ be the
measurable subset of $\Omega$ on which
$Y(t)\in K_{\epsilon}$.
We have
\begin{multline*}
\biggl(\expect\norm{ f_2(t,Y(t)) }^2\biggr)^{1/2}\\
\begin{aligned}
\leq&\min_{1\leq i\leq m}\biggl(\expect\bigl(
      \un{\Omega_{\epsilon,t}}  \norm{ f_2(t,Y(t))-f(t,y_i) }^2\bigr)\biggr)^{1/2}
     +\max_{1\leq i\leq m}\norm{f_2(t,y_i)}\\
 &+\biggl(\expect\bigl(\un{\Omega_{\epsilon,t}^c}\norm{f_2(t,Y(t))}^2\bigr)\biggr)^{1/2}\\
\leq&\epsilon+\max_{1\leq i\leq m}\norm{f_2(t,y_i)}
    +\bigl(\prob\CCO{\Omega_{\epsilon,t}^c}\bigr)^{1/q}
            \biggl(\expect\norm{f_2(t,Y(t))}^{p}\biggr)^{2/p}\\
\leq&\epsilon+\max_{1\leq i\leq m}\norm{f_2(t,y_i)} +\epsilon^{1/q}\Mp.
\end{aligned}
\end{multline*}
A similar result holds for $\expect\norm{ g_2(t,Y(t))}^2$.
Let us denote
\begin{align*}
\fg(\epsilon)&=\epsilon+\epsilon^{1/q}\Mp,\\
\fm_\epsilon(t)&=\max_{1\leq i\leq m}\norm{f_2(s,y_i)},\\
\gm_\epsilon(t)&=\max_{1\leq i\leq m}\norm{g_2(s,y_i)}.
\end{align*}
Thanks to the ergodicity of $f_2$ and $g_2$, the functions
$\fm_\epsilon$ and $\gm_\epsilon$ are in $\ER(\R,\R,\mu)$.
We have
\begin{multline*}
\begin{aligned}
\frac{1}{\mu([-r,r])}\int_{[-r,r]}
   \biggl(\expect\biggl\Vert&\int_{-\infty}^t S(t-s)f_2(s,Y(s))\,ds\\
         +&\int_{-\infty}^t S(t-s)g_2(s,Y(s))\,dW(s)\biggr\Vert^2\biggr)^{1/2} d\mu(t)
\end{aligned}\\
\begin{aligned}
\leq&\frac{1}{\mu([-r,r])}\int_{[-r,r]}
\biggl(   \CCO{\int_{-\infty}^t e^{-2\delta (t-s)}\expect\norm{f_2(s,Y(s))}^2 ds
        }^{1/2} \\
    &+(\trace Q)^{1/2}
   \CCO{\int_{-\infty}^t
         e^{ -2\delta (t-s)}\expect\norm{g_2(s,Y(s))}^2\,ds}^{1/2} \biggr)d\mu(t)\\
\leq&\frac{1 }{ \mu([-r,r]) }  \int_{[-r,r]}
   \biggl(\int_{-\infty}^t e^{-2\delta (t-s)}
            \CCO{\fm_\epsilon(s)+\fg(\epsilon)}^2\, ds
        \biggr)^{1/2} d\mu(t)\\
  &+(\trace Q)^{1/2}\frac{1 }{ \mu([-r,r]) }  \int_{[-r,r]}
   \biggl(\int_{-\infty}^t e^{-2\delta (t-s)}
            \CCO{\gm_\epsilon(s)+\fg(\epsilon)}^2\, ds
        \biggr)^{1/2} d\mu(t)\\
\leq &\frac{1}{ \mu([-r,r]) } \int_{[-r,r]}
   \biggl(\int_{-\infty}^t e^{-2\delta (t-s)}
            \fm_\epsilon^2(s)\, ds
        \biggr)^{1/2} d\mu(t)\\
  &+(\trace Q)^{1/2}\frac{1 }{ \mu([-r,r]) }  \int_{[-r,r]}
   \biggl(\int_{-\infty}^t e^{-2\delta (t-s)}
            \gm_\epsilon^2(s)\, ds
        \biggr)^{1/2} d\mu(t)\\
 &+\frac{1+(\trace Q)^{1/2} }{ \mu([-r,r]) } \fg(\epsilon)
   \int_{[-r,r]}\biggl(\int_{-\infty}^t e^{-2\delta (t-s)}\, ds\biggr)^{1/2}d\mu(t).
\end{aligned}
\end{multline*}
In the right hand side of the last inequality,
the last term is arbitrarily small and both other terms converge to
0 when $r$ goes to $+\infty$,
thanks to Lemma \ref{lem:convolut}.
We have thus proved that
$Z_{n+1}$ is in
$\ER\bigl(\R,\ellp{2}(\prob,\h_2),\mu\bigr)$.
We deduce by induction that the sequence  $(Z_n)$
lies in
$\ER\bigl(\R,\ellp{2}(\prob,\h_2),\mu\bigr)$.

Now, the sequence $(X_n)$ converges to $X$ in
$\CUB\bigl(\R,\ellp{2}(\prob, \h_2)\bigr)$,
thus
$(Z_n)$ converges to $Z:=X-Y$ in
$\CUB\bigl(\R,\ellp{2}(\prob, \h_2)\bigr)$.
Let $\epsilon>0$, and let $n$ such that
$$\sup_{t\in\R}\bigl(\expect\norm{Z(t)-Z_n(t)}^2\bigr)^{1/2}\leq\epsilon.$$
We have
\begin{multline*}
\frac{1}{\mu([-r,r])}\int_{[-r,r]}\bigl(\expect\norm{Z(t)}^2\bigr)^{1/2}d\mu(t)
\\
\begin{aligned}
\leq&\frac{1}{\mu([-r,r])}
      \int_{[-r,r]}\bigl(\expect\norm{Z(t)-Z_n(t)}^2\bigr)^{1/2}d\mu(t)\\
    &+\frac{1}{\mu([-r,r])}
      \int_{[-r,r]}\bigl(\expect\norm{Z_n(t)}^2\bigr)^{1/2}d\mu(t)\\
\leq& {\epsilon}
     +\frac{1}{\mu([-r,r])}
      \int_{[-r,r]}\bigl(\expect\norm{Z_n(t)}^2\bigr)^{1/2}d\mu(t).
\end{aligned}
\end{multline*}
As $\epsilon$ is arbitrary,
this proves that
$Z\in\ER\bigl(\R,\ellp{2}(\prob,\h_2),\mu\bigr)$.
\finproof
%%%%%%%%%

\begin{remark}
We did not use in the proof of Theorem \ref{theo:main} the hypothesis that $f_2\in \ERUb(\R\times \h_2,
\h_2,\mu)$ and $g_2\in \ERUb(\R\times \h_2, L(\h_1, \h_2),\mu)$. Actually, we needed only to assume that, for
each $x\in\h_2$, $f_2(.,x)\in \ER(\R, \h_2,\mu)$ and $g_2(.,x)\in \ER(\R, L(\h_1, \h_2),\mu)$. By Remark
\ref{rem:ERUc} and the Lipschitz condition, this is equivalent to assume that $f_2\in \ERUc(\R\times \h_2,
\h_2,\mu)$ and $g_2\in \ERUc(\R\times \h_2, L(\h_1, \h_2),\mu)$.
\end{remark}

\begin{theorem}\label{theo: main-ap}
{\bf(Weighted pseudo almost periodic solution of an equation with
  weighted pseudo almost periodic coefficients)}
Assume the same hypothesis as in Theorem \ref{theo:main}, and that $f_1$ and $g_1$ are almost periodic with
respect to the first variable, uniformly with respect to the second variable in bounded sets. Then (with obvious
definitions) the process $Y$ of Theorem \ref{theo:main} is almost
periodic in 2-distribution, 
thus the process $X$ is $\mu$-pseudo almost periodic in $2$-distribution.
\end{theorem}
\proof
The proof is exactly the same as that of Theorem \ref{theo:main},
replacing Theorem \ref{theo:AA} by \cite[Theorem
3.1]{KMRF12averaging}.
\finproof
%%%%%%%%%

\paragraph{Acknowledgements} We thank an anonymous referee for
pointing out an error in a previous version.

%%%%%%%%%%%%%%%%%%%%%%%%%%%%%%%%%%%%%%%%%%%%%%%%%%%%
\def\polhk#1{\setbox0=\hbox{#1}{\ooalign{\hidewidth
  \lower1.5ex\hbox{`}\hidewidth\crcr\unhbox0}}} \def\cprime{$'$}
  \def\cprime{$'$}

%%%%%%%%%%%%%%%%%%%%%%%%%%%%%%%%%%%%%%%%%%%%%%%%%%%%
%\bibliographystyle{hplain}\bibliography{BMRF.bib}%

\begin{thebibliography}{10}

\bibitem{abbas14weyl}
Syed Abbas.
\newblock A note on {W}eyl pseudo almost automorphic functions and their
  properties.
\newblock {\em Math. Sci. (Springer)}, 6:Art. 29, 5, 2012.

\bibitem{Arnold-Tudor98}
Ludwig Arnold and Constantin Tudor.
\newblock Stationary and almost periodic solutions of almost periodic affine
  stochastic differential equations.
\newblock {\em Stochastics Stochastics Rep.}, 64(3-4):177--193, 1998.

\bibitem{bedouhene-mellah-prf2012}
Fazia Bedouhene, Omar Mellah, and Paul Raynaud~de Fitte.
\newblock Bochner-almost periodicity for stochastic processes.
\newblock {\em Stoch. Anal. Appl.}, 30(2):322--342, 2012.

\bibitem{bezandry-diagana07existence}
Paul~H. Bezandry and Toka Diagana.
\newblock Existence of almost periodic solutions to some stochastic
  differential equations.
\newblock {\em Appl. Anal.}, 86(7):819--827, 2007.

\bibitem{Billingsley99}
Patrick Billingsley.
\newblock {\em Convergence of probability measures}.
\newblock Wiley Series in Probability and Statistics: Probability and
  Statistics. John Wiley \& Sons Inc., New York, second edition, 1999.
\newblock A Wiley-Interscience Publication.

\bibitem{blot-mophou-nguerekata-pennequin2009}
J.~Blot, G.~M. Mophou, G.~M. N'Gu{\'e}r{\'e}kata, and D.~Pennequin.
\newblock Weighted pseudo almost automorphic functions and applications to
  abstract differential equations.
\newblock {\em Nonlinear Anal.}, 71(3-4):903--909, 2009.

\bibitem{blot-cieutat_ezzinbi2012}
Jo{\"e}l Blot, Philippe Cieutat, and Khalil Ezzinbi.
\newblock Measure theory and pseudo almost automorphic functions: new
  developments and applications.
\newblock {\em Nonlinear Anal.}, 75(4):2426--2447, 2012.

\bibitem{blot-cieutat-ezzinbi2013}
Jo{\"e}l Blot, Philippe Cieutat, and Khalil Ezzinbi.
\newblock New approach for weighted pseudo-almost periodic functions under the
  light of measure theory, basic results and applications.
\newblock {\em Appl. Anal.}, 92(3):493--526, 2013.

\bibitem{blot-al09superposition}
Jo{\"e}l Blot, Philippe Cieutat, Gaston~M. N'Gu{\'e}r{\'e}kata, and Denis
  Pennequin.
\newblock Superposition operators between various almost periodic function
  spaces and applications.
\newblock {\em Commun. Math. Anal.}, 6(1):42--70, 2009.

\bibitem{bochner33abstrakte}
S.~Bochner.
\newblock Abstrakte {F}astperiodische {F}unktionen.
\newblock {\em Acta Math.}, 61(1):149--184, 1933.

\bibitem{bochner61uniform}
S.~Bochner.
\newblock Uniform convergence of monotone sequences of functions.
\newblock {\em Proc. Nat. Acad. Sci. U.S.A.}, 47:582--585, 1961.

\bibitem{bochner62new_approach}
S.~Bochner.
\newblock A new approach to almost periodicity.
\newblock {\em Proc. Nat. Acad. Sci. U.S.A.}, 48:2039--2043, 1962.

\bibitem{bochner55betti}
Salomon Bochner.
\newblock Curvature and {B}etti numbers in real and complex vector bundles.
\newblock {\em Univ. e Politec. Torino. Rend. Sem. Mat.}, 15:225--253,
  1955--56.

\bibitem{Bourbaki71-I-IV}
N.~Bourbaki.
\newblock {\em \'{E}l\'ements de math\'ematique. {T}opologie g\'en\'erale.
  {C}hapitres 1 \`a 4}.
\newblock Hermann, Paris, 1971.

\bibitem{casarino00}
Valentina Casarino.
\newblock Almost automorphic groups and semigroups.
\newblock {\em Rend. Accad. Naz. Sci. XL Mem. Mat. Appl. (5)}, 24:219--235,
  2000.

\bibitem{cenusa-sacuiu80}
Gh. Cenu{\c{s}}{\u{a}} and I.~S{\u{a}}cuiu.
\newblock Some properties of random functions almost periodic in probability.
\newblock {\em Rev. Roumaine Math. Pures Appl.}, 25(9):1317--1325, 1980.

\bibitem{DaPrato-Tudor95}
G.~Da~Prato and C.~Tudor.
\newblock Periodic and almost periodic solutions for semilinear stochastic
  equations.
\newblock {\em Stochastic Anal. Appl.}, 13(1):13--33, 1995.

\bibitem{diagana06weighted}
Toka Diagana.
\newblock Weighted pseudo almost periodic functions and applications.
\newblock {\em C. R. Math. Acad. Sci. Paris}, 343(10):643--646, 2006.

\bibitem{diagana07book}
Toka Diagana.
\newblock {\em Pseudo almost periodic functions in {B}anach spaces}.
\newblock Nova Science Publishers, Inc., New York, 2007.

\bibitem{diagana08weighted}
Toka Diagana.
\newblock Weighted pseudo-almost periodic solutions to some differential
  equations.
\newblock {\em Nonlinear Anal.}, 68(8):2250--2260, 2008.

\bibitem{diagana09stepanov}
Toka Diagana.
\newblock Existence of pseudo-almost automorphic solutions to some abstract
  differential equations with {${\bf S}^p$}-pseudo-almost automorphic
  coefficients.
\newblock {\em Nonlinear Anal.}, 70(11):3781--3790, 2009.

\bibitem{diagana12generalized_stepanov}
Toka Diagana.
\newblock Evolution equations in generalized {S}tepanov-like pseudo almost
  automorphic spaces.
\newblock {\em Electron. J. Differential Equations}, pages No. 49, 19, 2012.

\bibitem{ding-deng-nguerekata14}
Hui-Sheng Ding, Chao Deng, and Gaston~M. N'Gu{\'e}r{\'e}kata.
\newblock Almost automorphic random functions in probability.
\newblock {\em Abstr. Appl. Anal.}, pages Art. ID 243748, 6, 2014.

\bibitem{Engelking89}
Ryszard Engelking.
\newblock {\em General topology}, volume~6 of {\em Sigma Series in Pure
  Mathematics}.
\newblock Heldermann Verlag, Berlin, second edition, 1989.
\newblock Translated from the Polish by the author.

\bibitem{fink68aa-ap}
A.~M. Fink.
\newblock Almost automorphic and almost periodic solutions which minimize
  functionals.
\newblock {\em T\^ohoku Math. J. (2)}, 20:323--332, 1968.

\bibitem{fu2012}
Miaomiao Fu.
\newblock Almost automorphic solutions for nonautonomous stochastic
  differential equations.
\newblock {\em J. Math. Anal. Appl.}, 393(1):231--238, 2012.

\bibitem{fu-chen13}
Miaomiao Fu and Feng Chen.
\newblock Almost automorphic solutions for some stochastic differential
  equations.
\newblock {\em Nonlinear Anal.}, 80:66--75, 2013.

\bibitem{fu-liu2010square-mean-aa}
Miaomiao Fu and Zhenxin Liu.
\newblock Square-mean almost automorphic solutions for some stochastic
  differential equations.
\newblock {\em Proc. Amer. Math. Soc.}, 138(10):3689--3701, 2010.

\bibitem{halanay87}
A.~Halanay.
\newblock Periodic and almost periodic solutions to affine stochastic systems.
\newblock In {\em Proceedings of the {E}leventh {I}nternational {C}onference on
  {N}onlinear {O}scillations ({B}udapest, 1987)}, pages 94--101, Budapest,
  1987. J\'anos Bolyai Math. Soc.

\bibitem{halanay-tudor-morozan87}
A.~Halanay, C.~Tudor, and T.~Morozan.
\newblock Tracking almost periodic signals under white noise perturbations.
\newblock {\em Stochastics}, 21(4):287--301, 1987.

\bibitem{KMRF12averaging}
Mikhail {Kamenski}, Omar {Mellah}, and Paul {Raynaud de Fitte}.
\newblock {Weak Averaging of Semilinear Stochastic Differential Equations with
  Almost Periodic Coefficients}.
\newblock {\em J. Math. Anal. Appl.}, 427 (1) (2015) 336-364, 2015.

\bibitem{liang-xiao-zhang10decomposition}
Jin Liang, Ti-Jun Xiao, and Jun Zhang.
\newblock Decomposition of weighted pseudo-almost periodic functions.
\newblock {\em Nonlinear Anal.}, 73(10):3456--3461, 2010.

\bibitem{liang-zhang_xiao-jun08}
Jin Liang, Jun Zhang, and Ti-Jun Xiao.
\newblock Composition of pseudo almost automorphic and asymptotically almost
  automorphic functions.
\newblock {\em J. Math. Anal. Appl.}, 340(2):1493--1499, 2008.

\bibitem{lindgren2006stationary_processes}
Georg Lindgren.
\newblock Lectures on stationary stochastic processes.
\newblock Technical report, Lund University, Lund, Sweden, 2006.
\newblock A course for {PHD} students in Mathematical Statistics and other
  fields available at
  http://www.maths.lth.se/matstat/staff/georg/Publications/lecture2006.pdf.

\bibitem{liu-sun2014}
Zhenxin Liu and Kai Sun.
\newblock Almost automorphic solutions for stochastic differential equations
  driven by {L}\'evy noise.
\newblock {\em J. Funct. Anal.}, 266(3):1115--1149, 2014.

\bibitem{lizama-nguerekata10}
Carlos Lizama and Gaston~M. N'Gu{\'e}r{\'e}kata.
\newblock Bounded mild solutions for semilinear integro differential equations
  in {B}anach spaces.
\newblock {\em Integral Equations Operator Theory}, 68(2):207--227, 2010.

\bibitem{omarPHD}
Omar Mellah.
\newblock {\em \'Etude des solutions presque p\'eriodiques des \'equations
  diff\'erentielles stochastiques}.
\newblock PhD thesis, University Mouloud Mammeri, Tizi Ouzou, Algeria, and
  University of Rouen, France, 2012.

\bibitem{MRF13}
Omar Mellah and Paul Raynaud~de Fitte.
\newblock Counterexamples to mean square almost periodicity of the solutions of
  some {SDE}s with almost periodic coefficients.
\newblock {\em Electron. J. Differential Equations}, pages No. 91, 7, 2013.

\bibitem{Morozan-Tudor89}
T.~Morozan and C.~Tudor.
\newblock Almost periodic solutions of affine {I}t\^o equations.
\newblock {\em Stochastic Anal. Appl.}, 7(4):451--474, 1989.

\bibitem{nguerekata01book}
Gaston~M. N'Guerekata.
\newblock {\em Almost automorphic and almost periodic functions in abstract
  spaces}.
\newblock Kluwer Academic/Plenum Publishers, New York, 2001.

\bibitem{nguerekata05comments}
Gaston~M. N'gu{\'e}r{\'e}kata.
\newblock Comments on almost automorphic and almost periodic functions in
  {B}anach spaces.
\newblock {\em Far East J. Math. Sci. (FJMS)}, 17(3):337--344, 2005.

\bibitem{nguerekata05book}
Gaston~M. N'Gu{\'e}r{\'e}kata.
\newblock {\em Topics in almost automorphy}.
\newblock Springer-Verlag, New York, 2005.

\bibitem{nguerekata-pankov08stepanov}
Gaston~M. N'Gu{\'e}r{\'e}kata and Alexander Pankov.
\newblock Stepanov-like almost automorphic functions and monotone evolution
  equations.
\newblock {\em Nonlinear Anal.}, 68(9):2658--2667, 2008.

\bibitem{Onicescu-Istratescu75}
Octav Onicescu and Vasile~I. Istr{\u{a}}tescu.
\newblock Approximation theorems for random functions.
\newblock {\em Rend. Mat. (6)}, 8:65--81, 1975.
\newblock Collection of articles dedicated to Mauro Picone on the occasion of
  his ninetieth birthday.

\bibitem{Precupanu82}
A.~M. Precupanu.
\newblock On the almost periodic functions in probability.
\newblock {\em Rend. Mat. (7)}, 2(3):613--626, 1982.

\bibitem{shen-yi98semiflows}
Wenxian Shen and Yingfei Yi.
\newblock Almost automorphic and almost periodic dynamics in skew-product
  semiflows.
\newblock {\em Mem. Amer. Math. Soc.}, 136(647):x+93, 1998.

\bibitem{Slutsky38}
E.~Slutsky.
\newblock {Sur les fonctions al\'eatoires presque p\'eriodiques et sur la
  d\'ecomposition des fonctions al\'eatoires stationnaires en composantes.}
\newblock {Actual. sci. industr. 738, 33-55. (Conf\'er. internat. Sci. math.
  Univ. Gen\`eve. Th\'eorie des probabilit\'es. V: Les fonctions al\'eatoires.)
  (1938).}, 1938.

\bibitem{Tudor92affine}
C.~Tudor.
\newblock Almost periodic solutions of affine stochastic evolution equations.
\newblock {\em Stochastics Stochastics Rep.}, 38(4):251--266, 1992.

\bibitem{Tudor95ap_processes}
C.~Tudor.
\newblock Almost periodic stochastic processes.
\newblock In {\em Qualitative problems for differential equations and control
  theory}, pages 289--300. World Sci. Publ., River Edge, NJ, 1995.

\bibitem{Tudor-Tudor99}
C.~A. Tudor and M.~Tudor.
\newblock Pseudo almost periodic solutions of some stochastic differential
  equations.
\newblock {\em Math. Rep. (Bucur.)}, 1(51)(2):305--314, 1999.

\bibitem{veech65aa}
W.~A. Veech.
\newblock Almost automorphic functions on groups.
\newblock {\em Amer. J. Math.}, 87:719--751, 1965.

\bibitem{Whitt70}
Ward Whitt.
\newblock Weak convergence of probability measures on the function space
  {$C[0,\,\infty )$}.
\newblock {\em Ann. Math. Statist.}, 41:939--944, 1970.

\bibitem{xia-fan12}
Zhinan Xia and Meng Fan.
\newblock Weighted {S}tepanov-like pseudo almost automorphy and applications.
\newblock {\em Nonlinear Anal.}, 75(4):2378--2397, 2012.

\bibitem{zhang94integr}
Chuan~Yi Zhang.
\newblock Integration of vector-valued pseudo-almost periodic functions.
\newblock {\em Proc. Amer. Math. Soc.}, 121(1):167--174, 1994.

\bibitem{zhang94}
Chuan~Yi Zhang.
\newblock Pseudo-almost-periodic solutions of some differential equations.
\newblock {\em J. Math. Anal. Appl.}, 181(1):62--76, 1994.

\bibitem{zhang95}
Chuan~Yi Zhang.
\newblock Pseudo almost periodic solutions of some differential equations.
  {II}.
\newblock {\em J. Math. Anal. Appl.}, 192(2):543--561, 1995.

\bibitem{zhang01book}
Chuanyi Zhang.
\newblock {\em Almost periodic type functions and ergodicity}.
\newblock Science Press Beijing, Beijing; Kluwer Academic Publishers,
  Dordrecht, 2003.

\bibitem{zhang-chang-nguerekata12}
Rui Zhang, Yong-Kui Chang, and G.~M. N'Gu{\'e}r{\'e}kata.
\newblock New composition theorems of {S}tepanov-like weighted pseudo almost
  automorphic functions and applications to nonautonomous evolution equations.
\newblock {\em Nonlinear Anal. Real World Appl.}, 13(6):2866--2879, 2012.

\end{thebibliography}
%%%%%%%%%%%%%%%%%%%%%%%%%%%%%%%%%%%%%%%%%%%%%%%%%%%%

\end{document}